\documentclass[12pt]{amsart}

\usepackage[T1]{fontenc}
\usepackage[english]{babel}
\usepackage{amsmath,amssymb}
\usepackage{amscd}

\newtheorem{prop}{Theorem}[section]
\newtheorem{lemma}[prop]{Lemma}
\newtheorem{corollary}[prop]{Corollary}
\newtheorem{definition}[prop]{Definition}

\newtheorem{remark}[prop]{Remark}

\newcommand{\nn}{\mathbb{N}}
\newcommand{\zz}{\mathbb{Z}}
\newcommand{\rr}{\mathbb{R}}
\newcommand{\cc}{\mathbb{C}}
\newcommand{\inner}[1]{\langle #1 \rangle}
\newcommand{\Akom}{A\otimes\mathcal{K}}
\newcommand{\Bkom}{B\otimes\mathcal{K}}
\newcommand{\Kom}[1]{#1\otimes\mathcal{K}}
\newcommand{\kom}{\mathcal{K}}
\newcommand{\ope}{\mathcal{L}}
\newcommand{\mult}{\mathcal{M}}
\newcommand{\calc}{\mathcal{C}}
\newcommand{\diff}{\mathcal{D}}

\newcommand{\alg}{\mathcal{A}}

\newcommand{\grad}{\hat\otimes}
\newcommand{\gradodd}{\hat\otimes\varepsilon}
\newcommand{\spec}{\textrm{Sp}}
\newcommand{\modu}{\mathbb{E}}

\newcommand{\conn}{\quad\textrm{and}\quad}

\newcommand{\fin}{\mathcal{F}}
\newcommand{\diag}{\textrm{diag}}
\newcommand{\D}{\mathcal{D}}
\newcommand{\A}{\mathcal{A}}
\begin{document}

 \title{$KK$-Theory and Spectral Flow in von Neumann Algebras}

\author{J. Kaad, R. Nest, A. Rennie}
\thanks{email:
    \texttt{jenskaad@hotmail.com}, \texttt{rnest@math.ku.dk},
    \texttt{rennie@math.ku.dk}} 
\maketitle
\vspace{-10pt}
\centerline{ Institute for Mathematical
  Sciences, University of Copenhagen}
\centerline{ Universitetsparken 5, DK-2100 Copenhagen, Denmark}

\bigskip

\centerline{\textbf{Abstract}} 
We present a definition of spectral flow for any norm closed
ideal $J$ in any von Neumann algebra $N$. 
Given a path of selfadjoint operators in
$N$ which are invertible in $N/J$, the spectral flow produces a class
in $K_0(J)$. 

Given a semifinite spectral triple $(\alg,H,\D)$ relative to
$(N,\tau)$, 
we construct a class $[\diff]\in KK^1(A,\kom(N))$. 
For a unitary $u\in \alg$, the von Neumann spectral flow between
$\diff$  
and $u^*\diff u$ is equal to the Kasparov product
$[u]\grad_A[\diff]$, and is simply related to the numerical spectral
flow, and a  refined  $C^*$-spectral flow.
\vspace{-8pt}
\tableofcontents
\newpage
\section{Introduction}

The theory of analytic spectral flow has received a great deal of
attention in recent years, with significant progress being made by
many authors, \cite{BCPRSW, BLP,CP1,CP2,CPRS2,Ph,Ph1,PR,W}. The
article \cite{W} contains a much more detailed review of other aspects
of spectral flow.

Here we take a slightly different tack, replacing 
numerical measures of spectral flow by $K$-theory valued measures, as
in \cite{pr,W}. The
advantages of this approach are the great generality in which it can
be defined, and its compatibility with the various numerical notions.

This compatability yields constraints on the possible values of
spectral flow, which, for example, in the semifinite setting of
\cite{Ph,Ph1}, is a priori any real number. Our description of
spectral flow allows one to factor through a $K$-theory group, and so
constrain the possible values of the spectral flow. The more refined
we can be about the target $K$-theory group, the more refined our
information.

We define von Neumann spectral flow for any norm closed
ideal $J$ in any von Neumann algebra $N$. Given 
a path of selfadjoint operators in
$N$ which are invertible in $N/J$, we obtain a class in $K_0(J)$.
In order to be able to work in such a general context, we need to
develop a $K_0(J)$-valued index theory for any such pair $N,J$. Such
an index theory is developed in Section \ref{index-theory}, and then
in Section \ref{vN-spec-flow} we define and study the von Neumann
spectral flow. 
We then follow the approach of \cite{Ph,Ph1}, defining spectral flow in
terms of relative indices of projections.

A closely related idea which we introduce is a von Neumann spectral
triple, modelled on the definition of semifinite spectral triples, but
valid for any von Neumann algebra $N$ and ideal $J$. 
We show that such a triple defines
a Kasparov class, and relate the spectral flow to Kasparov products.

In particular, every semifinite spectral triple represents a
$KK$-class, just as ordinary spectral triples represent $K$-homology
classes. This extends the observed relation in \cite{pr,prs}.

In Section 5 we discuss the consequences of refining our target
$K$-theory group to $K_0({\mathcal B})$, where ${\mathcal B}\subset J$
is a $\sigma$-unital subalgebra. We show that this can always be done
for 
a von Neumann spectral triple, and so we can define a $C^*$-spectral
flow. We relate this spectral flow to our previously defined von
Neumann spectral flow.

Section 6 relates both von Neumann and $C^*$ spectral  flow for a
semifinite spectral triple to the numerical spectral flow obtained
from a trace. 

The Appendix summarises some results from $KK$-theory that we require,
and proves an explicit formula for certain odd pairings in
$KK$-theory, which plays a key role throughout the paper.

{\bf Acknowledgements} It is a pleasure to thank Alan Carey and 
John Phillips for many helpful conversations about spectral flow.

\section{K-theory-valued von Neumann Index Theory}\label{index-theory}

Throughout this section, we let 
$N$ be a von-Neumann algebra acting on a Hilbert space $H$ and let
$J$ be a norm closed ideal in $N$. Let $\pi : N\to N/J$ be the
quotient map.

In all the following, 
we will distinguish between the kernel of an operator $S: H\to H$
called $\ker(S)$ and the projection onto the kernel
called $N(S)\in\ope(H)$. Likewise we have the image of $S$,
$\mbox{Im}(S)$ and
the projection onto the norm closure of $\mbox{Im}(S)$, denoted 
$R(S)\in\ope(H)$. If $S$ is in 
$N$ then $N(S)$ and $R(S)$ are in $N$ also.

For any two projections $p,q\in \ope(H)$ we denote the projection onto
$\mbox{Im}(p)\cap \mbox{Im}(q)$  by $p\cap
q\in\ope(H)$. If $p,q\in N$ and $S\in N'$, then $Sp=pS$ and $Sq=qS$
thus $Sp\cap q=p\cap q S$. It follows easily that $p\cap q\in N''=N$.


We recall some facts about the polar decomposition of an operator.
Let $S\in N$. The partial isometry $u\in N$ from the polar
decomposition of 
$S$ is called the \emph{phase} of $S$. The phase of $S$ has the
following 
properties
\[
\begin{array}{cc}
u|S|=S & S^*=u^*|S^*|\\
uu^*=R(u)=R(S) & u^*u=R(u^*)=R(S^*)\\
1-uu^*=N(u^*)=N(S^*) & 1-u^*u=N(u)=N(S)
\end{array}
\]
See \cite[Appendice III]{Dixmier} for more details. Since $K$-theory
is well-defined for non-separable $C^*$-algebras, we can ask what the
generalised index map in $K$-theory gives us for an invertible in the
`Calkin algebra' $N/J$.


\begin{lemma}\label{connection}
Let $[\pi(S)]\in K_1(N/J)$ be a class in $K_1(N/J)$
represented by the unitary $\pi(S)$, where $S\in M_n(N)$ for some
$n\in \nn$. Then
\[
\partial[\pi(S)]=[N(S)]-[N(S^*)]\in K_0(J),
\]
where $\partial: K_1(N/J)\to K_0(J)$ is the boundary map in $K$-theory.
See for instance \cite[Definition $8.3.1$]{Black}.
\end{lemma}
\begin{proof}
The algebra $M_n(N)=M_n(\mathbb{C})\otimes N$ 
is a von-Neumann algebra acting on the Hilbert space
$\oplus_{i=1}^n H$ so $S$ can be polar decomposed
in $M_n(N)$. Let $u\in M_n(N)$ be the phase of $S$.
Now, $u$ is a lift of $\pi(S)$ since
\[
\pi(S)=\pi(u|S|)=\pi(u)\pi(S^*S)^{1/2}=\pi(u)
\]
And we conclude from the definition of the boundary map,
\cite{Black,Higson,Roerdam},  that
\[
\partial[\pi(S)]=[1-u^*u]-[1-uu^*]=[N(S)]-[N(S^*)]
\]
as claimed.
\end{proof}

The generic situation where the index of an operator $S$ is relevant for
applications is when $S:H_1\to H_2$. Even to define `odd' index
pairings one requires such operators. Thus one must consider operators
not in a von Neumann algebra $N$, but in a skew-corner $qNp$ where
$p,q\in N$ are projections. This situation was first considered in
\cite{CPRS3} for semifinite von Neumann algebras. The following
definition generalises the semifinite notion of Fredholm.

\begin{definition}
Let $p,q$ be projections in $N$. Then $S\in qNp$ is a
$(q$-$p)$-\emph{Fredholm operator} if there exists a $T,R\in pNq$ such that
\[
\pi(TS)=\pi(p)\conn \pi(SR)=\pi(q)
\]
Since $\pi(T)=\pi(TSR)=\pi(R)$,  we can choose $R=T$. The operator $T$ 
is called a \emph{parametrix} for $S$. 
\end{definition}

\begin{remark}\label{commutproj}
Suppose we have an operator $S\in qNp$. Then 
\[
N(S)\cap p = N(S)p \conn N(S^*)\cap q=N(S^*)q
\]
This follows immediately since
\begin{equation*}
(1-p)N(S)=(1-p)=N(S)(1-p)\Rightarrow pN(S)=N(S)p
\end{equation*}
so $N(S)\cap p=N(S)p$.
Similar comments apply to the projections $N(S^*)$ and $q$.
\end{remark}

\begin{lemma}\label{polarfred}
Let $S\in qNp$ and let $u\in N$ be the phase of $S$.
Then $u\in qNp$ and we have the identities
\[
\begin{array}{c}
p-u^*u=N(S)-(1-p)=N(S)\cap p \\
q-uu^*=N(S^*)-(1-q)=N(S^*)\cap q.
\end{array}
\]
Furthermore if $S\in qNp$ is a $(q$-$p)$-Fredholm operator then $\pi(u^*u)=\pi(p)$ and
$\pi(uu^*)=\pi(q)$. In particular $u$ is $(q$-$p)$-Fredholm and 
$N(S)\cap p, N(S^*)\cap q\in J$
\end{lemma}
\begin{proof}
First, $u$ is in $qNp$ since $(1-p)H\subseteq\textrm{Ker}(S)=\textrm{Ker}(u)$
and $\textrm{Im}(u)=\overline{\textrm{Im}(S)}\subseteq qH$. Next,
$(1-p)N(S)=(1-p)$ so $N(S)-(1-p)=N(S)-N(S)(1-p)=N(S)p=N(S)\cap p$ by
Remark \ref{commutproj}.
The statement concerning $N(S^*)$ and $q$ is proved in the same way.

Now, suppose that $S\in qNp$ is a $(q$-$p)$-Fredholm operator with parametrix 
$T\in pNq$. Then $S^*S\in pNp$ is a $(p$-$p)$-Fredholm operator with parametrix $TT^*\in pNp$.
This means that $\pi(S^*S)$ is invertible in the $C^*$-algebra $\pi(p)N/J\pi(p)$. 
Similarly $\pi(SS^*)$ is invertible in the $C^*$-algebra $\pi(q)N/J\pi(q)$
Clearly, then the phase $u\in qNp$ of $S\in qNp$ is a lift of 
$\pi(S)\pi(S^*S)^{-1/2}\in \pi(q)N/J\pi(p)$. This allows us to deduce the identities
\[
\begin{array}{cc}
\pi(u^*u)=\pi(S^*S)^{-1/2}\pi(S^*S)\pi(S^*S)^{-1/2}=\pi(p) & \textrm{and} \\
\pi(uu^*)=\pi(S)\pi(S^*S)^{-1}\pi(S^*)=\pi(q)
\end{array}
\]
as desired.
\end{proof}

The result allows us to make the following definition.
\begin{definition}\label{indexdef}
Let $S\in qNp$ be a $(q$-$p)$-Fredholm operator.
We define the $(q$-$p)$-\emph{index} of $S$ as the class
\[
\textrm{\emph{Ind}}_{(q\textrm{\emph{-}p})}(S)=[N(S)\cap p]-[N(S^*)\cap q]
\]
in $K_0(J)$.
\end{definition}

Let $S\in qNp$ be a $(q$-$p)$-Fredholm operator and let $u\in qNp$ be the 
phase of $S$. The triple $(p,q,u)$ is a relative $K$-cycle and thus defines 
the class $[S]:=[p,q,u]\in K_0(N,N/J)$ in relative $K$-theory. The relative 
$K$-theory $K_0(N,N/J)$ is related to the $K$-theory of the ideal $K_0(J)$ 
through the excision map 
\[
\textrm{Exc} : K_0(J)\to K_0(N,N/J)
\]
as defined in \cite[Definition $4.3.7$]{Higson}. 
The excision map is an isomorphism, \cite[Theorem $4.3.8$]{Higson}.
In the next theorem we will see that the $(q$-$p)$-index of $S$ is simply the inverse 
of the excision map 
applied to the class $[S]\in K_0(N,N/J)$. Many properties of the $(q$-$p)$-index 
will thus follow immediately, and we will state the ones we need as corollaries.

\begin{prop}\label{indexc}
Let $S\in qNp$ be a $(q$-$p)$-Fredholm operator and let $u\in qNp$ be the 
phase of $S$. Then the identity
\[
\emph{Exc}^{-1}[S]=\emph{Ind}_{q\textrm{-}p}(S)
\]
is valid in $K_0(J)$
\end{prop}
\begin{proof}
We can express the class $[S]\in K_0(N,N/J)$ as a sum of classes
\[
[S]=[p,q,u]=[p-u^*u,q-uu^*,0]+[u^*u,uu^*,u]
\]
The relative $K$-cycle $(u^*u,uu^*,u)$ is degenerate so actually 
\[
[S]=[p-u^*u,q-uu^*,0]
\]
The projections $p-u^*u=N(S)\cap p$ and $q-uu^*=N(S^*)\cap q$ are in $J$ by 
Lemma \ref{polarfred}, so 
\[
\textrm{Exc}^{-1}[S]=[p-u^*u]-[q-uu^*]=\textrm{Ind}_{q\textrm{-}p}(S)
\]
as desired.
\end{proof}

\begin{corollary}\label{homoinvar}
Let $S_0\in qNp$ and $S_1\in qNp$ be $(q$-$p)$-Fredholm operators.
Suppose that there is a norm-continuous path of $(q$-$p)$-Fredholm operators
connecting $S_0$ and $S_1$. Then
\[
\textrm{\emph{Ind}}_{(q\textrm{\emph{-}p})}(S_0)=
\textrm{\emph{Ind}}_{(q\textrm{\emph{-}p})}(S_1)
\]
\end{corollary}
\begin{proof}
Let $t\mapsto S_t\in qNp$ be the norm-continuous path connecting $S_0$ and $S_1$.
The norm-continuous path $t\mapsto \pi(S_t)\pi(S_t^*S_t)^{-1/2}\in \pi(q)N/J\pi(p)$, where the 
inverse is in $\pi(p)N/J\pi(p)$, lifts to a 
path $t\mapsto v_t\in qN/Jp$ such that $(p,q,v_t)$ are relative $K$-cycles for all 
$t\in [0,1]$, \cite[Lemma $4.3.13$]{Higson}. If $u_0\in qNp$ and $u_1\in qNp$ are the phases of $S_0$ and 
$S_1$ respectively, then $\pi(u_0)=\pi(v_0)$ and $\pi(u_1)=\pi(v_1)$ so we have 
the identity 
\[
[S_0]=[p,q,u_0]=[p,q,v_0]=[p,q,v_1]=[p,q,u_1]=[S_1]
\]
in $K_0(N,N/J)$. It thus follows immediately by Theorem \ref{indexc} that
\[
\textrm{Ind}_{(q\textrm{-}p)}(S_0)=\textrm{Exc}^{-1}[S_0]=\textrm{Exc}^{-1}[S_1]=
\textrm{Ind}_{(q\textrm{-}p)}(S_1)
\]
as desired.
\end{proof}

\begin{corollary}\label{algebraic}
Let $S\in qNp$ be a $(q$-$p)$-Fredholm operator and let $T\in rNq$ be an
$(r$-$q)$-Fredholm operator. Then $TS$ is an $(r$-$p)$-Fredholm operator and
\[
\emph{Ind}_{(r\textrm{-}q)}(T)+
\emph{Ind}_{(q\textrm{-}p)}(S)=
\emph{Ind}_{(r\textrm{-}p)}(TS)
\]
\end{corollary}
\begin{proof}
Let $v\in rNq$, $u\in qNp$ and $w\in rNp$ be the phases of $T$, $S$ 
and $TS$ respectively. From the calculation
\[
\begin{split}
\pi(w)&=\pi(TS)\pi(S^*T^*TS)^{-1/2}\\
&=\pi(T)\pi(SS^*T^*T)^{-1/2}\pi(S)\\
&=\pi(T)\pi(T^*T)^{-1/2}\pi(S)\pi(S^*S)^{-1/2}\\
&=\pi(vu)
\end{split}
\]
we deduce the identity $[p,r,vu]=[p,r,w]$ in $K_0(N,N/J)$.

Summing the classes $[T]$ and $[S]$ in $K_0(N,N/J)$ we get
\[
[T]+[S]=[q,r,v]+[p,q,u]=[p,r,vu]=[p,r,w]=[TS]
\]
where the second equality follows from the relations in $K_0(N,N/J)$. 

This allows us to conclude that 
\[
\textrm{Ind}_{(r\textrm{-}q)}(T)+\textrm{Ind}_{(q\textrm{-}p)}(S)=
\textrm{Exc}^{-1}[T]+\textrm{Exc}^{-1}[S]=\textrm{Exc}^{-1}[TS]=
\textrm{Ind}_{(r\textrm{-}p)}(TS)
\]
as desired.
\end{proof}

Let $N$ be a semifinite von Neumann algebra equipped with a fixed normal, 
semifinite, faithful trace $\tau$. Let $\kom_N$ be the $\tau$-compact 
operators as defined in Definition \ref{taucompact}. All projections 
in $\kom_N$ have finite trace by Theorem \ref{fintrace}. Applying the 
homomorphism $\tau_* : K_0(\kom_N)\to\rr$ from Theorem \ref{homokom}
to the 
main theorems 
of this section we obtain some of the important results from Breuer-Fredholm 
theory, \cite{BCPRSW,B1,B2,CP1,CP2,CPRS2,CPRS3,Ph,Ph1,PR}.



\section{Von Neumann Spectral Flow}\label{vN-spec-flow}



\subsection{Basic Definitions and Properties}




Having set up the appropriate index theory for Fredholm operators in
skew-corners $pNq$, \cite{CPRS3}, we now analyse spectral flow. This
is associated with odd index pairings, and so self-adjoint
operators. Specialising our definition of $p$-$q$-Fredholm to the case
$p=q=1$ we have the following.

\begin{definition}
An operator $T\in N$ is said to be $J$\emph{-Fredholm} if $\pi(T)$
is invertible in $N/J$. The space of $J$-Fredholm operators is denoted by
$\mathcal{F}$. The space of selfadjoint $J$-Fredholm operators is
denoted by $\mathcal{F}_{sa}$
\end{definition}

Let $\chi : \rr\to \rr$ be the indicator function for the 
interval $[0,\infty)$ defined by
\[
\chi(t)=\left\{\begin{array}{cc}
  1 & \qquad t\in [0,\infty) \\
  0 & \qquad t\in (-\infty,0)
\end{array}\right.
\]

The following lemma was first proved in \cite[Lemma $4.3$]{BCPRSW} in
the semifinite context. In fact it makes sense and is true in the more
general context considered here. We quote the statement and proof for
completeness. 

\begin{lemma}\label{change}
Let $T\in \mathcal{F}_{sa}$ then
$\pi\big(\chi(T)\big)=\chi\big(\pi(T)\big)$.
\end{lemma}
\begin{proof}
Note that $\chi\big(\pi(T)\big)$ makes sense since
$0\notin\textrm{Sp}\big(\pi(T)\big)$ thus we can find an
$\varepsilon>0$ such that the interval $[-\varepsilon,\varepsilon]$ is 
included in the resolvent set of $\pi(T)$. Now define the continuous
functions $f_1 : \rr\to \rr$
\[
f_1(t)=\left\{\begin{matrix} 0 & \qquad t\in (-\infty,-\varepsilon] \\
\varepsilon^{-1}t + 1 & \qquad t\in
[-\varepsilon,0]\\
1 & \qquad t\in [0,\infty)
\end{matrix}\right.
\]
and $f_2:\rr\to\rr$ by
\[
f_2(t)=\left\{\begin{matrix} 0 & \qquad t\in (-\infty,0] \\
\varepsilon^{-1}t & \qquad t\in [0,\varepsilon]\\
1 & \qquad t\in [\varepsilon,\infty)
\end{matrix}\right.
\]
So $f_1=\chi=f_2$ on $\textrm{Sp}(\pi(T))$ while $f_1\geq\chi\geq
f_2$ on $\textrm{Sp}(T)$. Thus
\[
\chi\big(\pi(T)\big)=f_1\big(\pi(T)\big)= \pi\big(f_1(T)\big)\geq
\pi\big(\chi(T)\big)\geq\pi\big(f_2(T)\big)=f_2\big(\pi(T)\big)
=\chi\big(\pi(T)\big)
\]
yielding $\chi\big(\pi(T)\big)=\pi\big(\chi(T)\big)$ as desired.
\end{proof}


\begin{lemma}\label{contin}
Let $t\mapsto B_t$ be a norm continuous path in $\mathcal{F}_{sa}$.
Then $t\mapsto \chi\big(\pi(B_t)\big)$ is a norm continuous path in the
$C^*$-algebra $N/J$.
\end{lemma}
To prove this lemma we need a general result from the theory of
$C^*$-algebras. The result is probably well-known to the experts, but
as we could not find a reference, we include a proof.


\begin{lemma}\label{stjerne}
Let $A$ be a $C^*$-algebra and let $U$ be an open subset of $\rr$.
Denote by $A_{sa}$ the real subspace of selfadjoint elements with
the induced topology from $A$.
Then the set
\[
\{a\in A_{sa}\,|\,\textrm{\emph{Sp}}(a)\subseteq U\,\}
\]
is open in $A_{sa}$
\end{lemma}
\begin{proof}
Let $a\in A_{sa}$ with $\textrm{Sp}(a)\subseteq U$. The function
$\textrm{dist}(\cdot, U^c) : \cc\to [0,\infty[$ defined by
\[
\textrm{dist}(\lambda, U^c)=\inf\{|\lambda-\mu|\,|\,\mu\in U^c\}
\]
for all $\lambda\in\cc$ is continuous. It attains thus its minimum
on the compact set $\textrm{Sp}(a)$. Furthermore for
$\lambda\in\spec(a)$ we have $\textrm{dist}(\lambda, U^c)>0$ because
$\lambda\notin \overline{U^c}=U^c$, so the minimum is strictly
positive. Set
\[
\varepsilon = \frac{\textrm{dist}\big(\spec(a),U^c\big)}{2} =
\frac{\inf\{|\lambda-\mu|\,|\,\lambda\in\textrm{Sp}(a),\mu\in
U^c\}}{2}>0
\]
Now take $b\in A_{sa}$ with $\|b-a\|<\frac{\varepsilon}{2}$ and suppose
for contradiction that there exists a $\lambda\in \spec(b)$ with
\[
B_\varepsilon(\lambda)\cap \spec(a)=\emptyset
\]
Here $B_\varepsilon(\lambda)$ denotes the ball of radius
$\varepsilon>0$ and center $\lambda$. Let $\mu\in
B_{\varepsilon/4}(\lambda)$. Then $\mu\notin \spec(a)$ and
\[
\|(\mu-a)^{-1}\|^{-1}=
\sup\{|\mu-\alpha|^{-1}\,|\,\alpha\in\spec(a)\}^{-1}=
\textrm{dist}\big(\mu,\spec(a)\big)\geq\frac{3\varepsilon}{4}
\]
Furthermore
\[
\|(\lambda-b)-(\mu-a)\|\leq|\lambda-\mu|+\|a-b\|<
\frac{\varepsilon}{4}+\frac{\varepsilon}{2}\leq
\|(\mu-a)^{-1}\|^{-1}
\]
So $\lambda-b$ is actually invertible which is a contradiction,
see \cite[Proposition $17.3$]{Meise}.
Hence for $\lambda\in \spec(b)$ we cannot have
\[
B_\varepsilon(\lambda)\cap \spec(a)=\emptyset
\]
Because of the way the $\varepsilon$ was chosen we conclude that
$\spec(b)\subseteq U$. Thus $\spec(b)\subseteq U$ for any
$b\in A_{sa}$ with $\|b-a\|<\varepsilon/2$
\end{proof}


\begin{proof}
{\em of Lemma \ref{contin}.} Let $t_0\in [0,1]$. Choose an $\varepsilon
>0$ such that the interval $[-\varepsilon,\varepsilon]$ is included in the 
resolvent set of $\pi(B_{t_0})$. Now
\[
\spec\big(\pi(B_{t_0})\big)\subseteq
(-\infty,-\varepsilon)\cup(\varepsilon,\infty)
\]
By Lemma \ref{stjerne} and the continuity of $t\mapsto \pi(B_t)$ there is
a $\delta>0$ such that
\[
\spec\big(\pi(B_t)\big)\subseteq
(-\infty,-\varepsilon)\cup(\varepsilon,\infty)
\]
for all $t\in(t_0-\delta,t_0+\delta)\cap[0,1]$. So for all
$t\in(t_0-\delta,t_0+\delta)\cap[0,1]$ we have the identity
\[
\chi\big(\pi(B_t)\big)=f\big(\pi(B_t)\big)
\]
where $f$ is some fixed continuous function (for instance the
function $f_1$ from the proof of Lemma \ref{change}). But the function
\[
t\mapsto f\big(\pi(B_t)\big)
\]
is clearly continuous and the lemma is thereby proved.
\end{proof}

With these tools at hand we can now define spectral flow as a class in
$K_0(J)$.
\begin{definition}[Spectral flow]
Let $t\mapsto B_t$ be a norm continuous path in $\mathcal{F}_{sa}$.
By Lemma \ref{contin} the path
\[
t\mapsto\pi\big(\chi(B_t)\big)=\chi\big(\pi(B_t)\big)
\]
is norm continuous. Find a partition $0=t_0<t_1<\ldots<t_n=1$ such
that
\[
\|\pi\big(\chi(B_t)\big)-\pi\big(\chi(B_s)\big)\|<1/2\quad\textrm{for
all}\quad t,s\in[t_{i-1},t_i]
\]
Set $p_i=\chi(B_{t_i})$. We now define the \emph{spectral flow} of
the path $\{B_t\}$ to be
\[
\textrm{\emph{sf}}\{B_t\}=\sum_{i=1}^n\big[(1-p_i)\cap
p_{i-1}\big]-\big[(1-p_{i-1})\cap p_i\big]\in K_0(J)
\]
\end{definition}

This definition raises several questions which we will answer 
in the following lemmas.
\begin{enumerate}
\item Are the elements $p_ip_{i-1}\in p_iNp_{i-1}$ $(p_i$-$p_{i-1})$-Fredholm
operators for all $i\in \{1,\ldots,n\}$ ?
\item Is the spectral flow independent of the partition chosen ?
\item Is the spectral flow invariant under homotopies of the path
$\{B_t\}$ ?
\item Is the spectral flow of $\{B_t\}$ equal to the spectral flow of
$\{C_t\}$ if $B_t-C_t\in J$ for all $t\in [0,1]$ ?
\end{enumerate}


\begin{lemma}\label{Jcontained} Suppose that $p, q\in N$ are two projections
such that $\|\pi(p)-\pi(q)\|<1$. Then $qp\in qNp$ is a $(q$-$p)$-Fredholm operator.
Thus, by Lemma \ref{polarfred}, we have $(1-q)\cap p\in J$ and $(1-p)\cap q\in J$.
\end{lemma}
\begin{proof}
The inequality
\[
\|\pi(pqp)-\pi(p)\|\leq \|\pi(p)-\pi(q)\|<1
\]
shows that $\pi(pqp)$ is invertible in $\pi(pNp)$, so there is an operator
$T\in pNp$ such that $\pi(Tpqp)=\pi(p)$.
Likewise the inequality
\[
\|\pi(qpq)-\pi(q)\|\leq\|\pi(q)-\pi(q)\|<1
\]
shows that $\pi(qpq)$ is invertible in $\pi(qNq)$, so there is
an operator $R\in qNq$ such that $\pi(qpqR)=\pi(q)$.
It follows that $qp$ is a $(q$-$p)$-Fredholm operator.
\end{proof}


\begin{corollary}\label{beginning}
For a path $\{B_t\}$ in $\mathcal{F}_{sa}$ and a
partition $0=t_0<t_1<\ldots<t_n=1$ such that
\[
\|\pi\big(\chi(B_t)\big)-\pi\big(\chi(B_s)\big)\|<1/2
\qquad\textrm{for all }t,s\in[t_{i-1},t_i]
\]
for all $i\in \{1,\ldots,n\}$ we can express the spectral flow of the path as
the sum of $(p_i$-$p_{i-1})$-indices
\[
\textrm{\emph{sf}}\{B_t\}=\sum_{i=1}^n\textrm{\emph{Ind}}
_{(p_i\textrm{\emph{-}}p_{i-1})}(p_ip_{i-1})
\]
where $p_i=\chi(B_{t_i})$ for all $i\in \{0,\ldots,n\}$. Thus by Theorem
\ref{algebraic} we actually have
\[
\textrm{\emph{sf}}\{B_t\}
=\textrm{\emph{Ind}}_{(p_n\textrm{\emph{-}}p_0)}(p_n\ldots p_0)
=[N(p_n\ldots p_0)\cap p_0]-[N(p_0\ldots p_n)\cap p_n]
\]
\end{corollary}


\begin{lemma}\label{conclusion} Suppose that $p,q,r$ are three
projections in $N$ with
\[
\|\pi(p)-\pi(q)\|< 1/2\quad,\quad \|\pi(q)-\pi(r)\|< 1/2\conn \|\pi(r)-\pi(p)\|<1/2
\]
then
\[
\textrm{\emph{Ind}}_{(r\textrm{\emph{-}}q)}(rq)+
\textrm{\emph{Ind}}_{(q\textrm{\emph{-}}p)}(qp)=
\textrm{\emph{Ind}}_{(r\textrm{\emph{-}}p)}(rp)
\]
Thus the spectral flow is independent of the partition chosen -it
doesn't change if a finer one is chosen.
\end{lemma}
\begin{proof}
We want to prove that
\[
\textrm{Ind}_{(r\textrm{\emph{-}}q)}(rq)+
\textrm{Ind}_{(q\textrm{\emph{-}}p)}(qp)-
\textrm{Ind}_{(r\textrm{\emph{-}}p)}(rp)=0
\]
By Theorem \ref{algebraic} this amounts to show that
\[
\textrm{Ind}_{(r\textrm{\emph{-}}r)}(rqpr)=0
\]
Verify the inequality
\[
\begin{split}
\|\pi(rqpr)-\pi(r)\|&\leq\|\pi(qp)-\pi(r)\|\\
&\leq\|\pi(qp)-\pi(q)\|+\|\pi(q)-\pi(r)\|\\
&\leq \|\pi(p)-\pi(q)\|+\|\pi(q)-\pi(r)\|\\
&<1
\end{split}
\]
Let $t\in[0,1]$, then 
\[
\|\pi\big((1-t)rqpr+tr\big)-\pi(r)\|=
(1-t)\|\pi(rqpr)-\pi(r)\|<(1-t)
\]
thus $\pi\big((1-t)rqpr+tr\big)$ is invertible in $\pi(rNr)$ for all
$t\in [0,1]$. This means that the path $t\mapsto (1-t)rqpr+tr$ consists entirely of
$(r$-$r)$-Fredholm operators and it connects $rqpr$ with $r$. To finish the proof we simply
refer to Theorem \ref{homoinvar} which gives
\[
0=\textrm{Ind}_{(r\textrm{\emph{-}}r)}(r)=
\textrm{Ind}_{(r\textrm{\emph{-}}r)}(rqpr)
\]
as desired.
\end{proof}


\begin{lemma}\cite{BCPRSW,Ph1}\label{homotopy}
Let $\{B_t\}$ and $\{C_t\}$ be two paths of selfadjoint
$J$-Fredholm operators. Let $H : [0,1]\times [0,1]\to
\mathcal{F}_{sa}$ be a homotopy connecting $\{B_t\}$ and $\{C_t\}$
leaving the endpoints fixed. That is $H$ is norm-continuous with
$H(t,0)=B_t$, $H(t,1)=C_t$ for all $t\in [0,1]$ and $H(0,s)=B_0$,
$H(1,s)=B_1$ for all $s\in [0,1]$. In particular $B_0=C_0$ and
$B_1=C_1$. Then
$\textrm{\emph{sf}}\{B_t\}=\textrm{\emph{sf}}\{C_t\}$.
\end{lemma}
\begin{proof}
The map $\zeta : [0,1]\times[0,1]\to N/J$ defined by
\[
\zeta(t,s)=\pi\Big(\chi\big(H(t,s)\big)\Big)
\]
is continuous and thus uniformly continuous, so we can choose a grid
\[
0=t_0<t_1\ldots<t_n=1\quad,\quad 0=s_0<s_1\ldots<s_n=1
\]
of $[0,1]\times[0,1]$ such that for any $(t,s),(u,v)\in
[t_{i-1},t_i]\times[s_{j-1},s_j]$ we have
$\|\zeta(t,s)-\zeta(u,v)\|<\frac{1}{2}$ where $i,j\in\{1,\ldots,n\}$
are fixed.

Now look at the spectral flow along the borders of the squares. That
is, for $i,j\in\{1,\ldots,n\}$ there are eight paths of selfadjoint
$J$-Fredholm operators. For instance we have
\[
u\mapsto H\big((1-u)t_{i-1}+ut_i,s_{j-1}\big)
\]
as one of them. The spectral flow of this path will be denoted by 
\[
\textrm{sf}_H\big((t_{i-1},s_{j-1}),(t_i,s_{j-1})\big)
\]
Likewise for the spectral flow of the other paths. Applying
Lemma \ref{conclusion} and the definition of spectral flow gives
\[
\begin{split}
&\textrm{sf}_H\big((t_{i-1},s_{j-1}),(t_i,s_{j-1})\big)+ \textrm{sf}_H\big((t_i,s_{j-1}),(t_i,s_j)\big)\\
&\qquad+\textrm{sf}_H\big((t_i,s_j),(t_{i-1},s_j)\big)+
\textrm{sf}_H\big((t_{i-1},s_j),(t_{i-1},s_{j-1})\big)=0
\end{split}
\]
Furthermore
\[
sf_H\big((t_{i-1},s_{j-1}),(t_i,s_{j-1})\big)=-sf_H\big((t_i,s_{j-1}),(t_{i-1},s_{j-1})\big)
\]
And an easy combinatorial argument yields the result.
\end{proof}


\begin{remark}\label{zeroflow}
Suppose that $p,q\in N$ are two projections with $\|p-q\|<1$, then
\[
\emph{Ker}(p)\cap\emph{Im}(q)=0=\emph{Ker}(q)\cap\emph{Im}(p)
\]
so the $J$-index of the projections
\[
\emph{Ind}_{(p\emph{-}q)}(pq)=[(1-p)\cap q]-[(1-q)\cap p]=0
\]

To see this we start by deducing that $1-p+pqp$ is invertible in $N$
from the inequality
\[
\|p-pqp\|\leq\|p-q\|<1
\]
If now $x$ is in $\emph{Ker}(q)\cap\emph{Im}(p)$ we immediately have
\[
(1-p+pqp)x=0
\]
but $1-p+pqp$ was invertible so $x=0$. Therefore
$\emph{Ker}(q)\cap\emph{Im}(p)=0$.
To prove that $\emph{Ker}(p)\cap\emph{Im}(q)=0$
simply interchange $p$ and $q$.
\end{remark}


\begin{lemma}\label{pertubation}
Let $\{B_t\}$ and $\{C_t\}$ be two paths of self adjoint $J$-Fredholm
operators with $B_t-C_t\in J$ for all $t\in [0,1]$ and
\[
\textrm{\emph{Ind}}_{(p_0\emph{-}q_0)}(p_0q_0)=
\textrm{\emph{Ind}}_{(q_1\emph{-}p_1)}(q_1p_1)=0\label{blah}
\]
where $p_0=\chi(B_0)$, $p_1=\chi(B_1)$, 
$q_0=\chi(C_0)$ and $q_1=\chi(C_1)$.
Then $\textrm{\emph{sf}}\{B_t\}=\textrm{\emph{sf}}\{C_t\}$. The 
condition \eqref{blah} is true if for instance
\[
\|\chi(B_0)-\chi(C_0)\|<1\conn
\|\chi(C_1)-\chi(B_1)\|<1
\]
by Remark \ref{zeroflow}.
\end{lemma}
\begin{proof}
Choose a partition $0=t_0<t_1<\ldots<t_n=1$ such that
\[
\|\pi\big(\chi(B_t)\big)-\pi\big(\chi(B_s)\big)\|<\frac{1}{4}
\]
and
\[
\|\pi\big(\chi(C_t)\big)-\pi\big(\chi(C_s)\big)\|<\frac{1}{4}
\]
for all $t,s\in[t_{i-1},t_i]$, $i\in\{1,\ldots,n\}$.

Now join the elements $B_{t_i}$ and $C_{t_i}$ by a straight line for
each $i\in\{0,\ldots,n\}$ denoted by $(BC)_i$. The straight line
from $C_{t_i}$ to $B_{t_i}$ is denoted by $(CB)_i$.

Notice that the lines are paths of selfadjoint $J$-Fredholm
operators because
\[
\pi\big((1-t)B_{t_i}+tC_{t_i}\big)=\pi(B_{t_i})
\]
for all $t\in[0,1]$ and $i\in\{1,\ldots,n\}$.

Now, almost by definition, the spectral flow along the square
\[
\begin{CD}
C_{t_{i-1}}@<<< C_{t_i} \\
@V{(CB)_{i-1}}VV @A{(BC)_i}AA \\
B_{t_{i-1}}@>>> B_{t_i}
\end{CD}
\]
is zero. Since too the spectral flow along the lines $(BC)_0$ and
$(BC)_1$ is zero by assumption we can use the same combinatorial
argument as in the proof of Lemma \ref{homotopy} to reach the desired
conclusion, namely $sf\{B_t\}=sf\{C_t\}$.
\end{proof}



\subsection{Von Neumann 
Spectral Triples and Spectral Flow}



\begin{definition}\label{specdef}
A \emph{von Neumann spectral triple}  $(\alg, H, \diff)$ relative to
$(N,J)$ 
consists of a representation of the $*$-algebra $\alg$ in the von
Neumann algebra $N$ acting on the Hilbert space $H$, together with a
norm closed ideal $J$ and a self-adjoint operator $\D$ affiliated to
$N$ such that
\begin{enumerate}
\item $[\diff, a]$ is defined on $\textrm{Dom}(\diff)$ and extends to a bounded
operator on $H$ for all $a\in\alg$.
\item $a(\lambda-\diff)^{-1}\in J$ for all $\lambda\notin\rr$ and $a\in \alg$.
\end{enumerate}
The $J$-spectral triple is said to be \emph{unital} if the unit of
$N$ is in $\alg$.
\end{definition}


If $(\alg, H, \diff)$ is a unital $J$-spectral triple, we use the spectral
theorem to define the bounded
operator in $N$
\[
F_\diff:=\diff(1+\diff^2)^{-1/2}.
\]

Let $t\mapsto A_t$
be a path of selfadjoint operators in $N$. 
We claim the path
\[
t\mapsto \diff_t:=\diff + A_t
\]
is a continuous path of unbounded selfadjoint $J$-Fredholm
operators in the sense that the path
\[
t\mapsto F_{\diff_t}=\diff_t(1+\diff_t^2)^{-\frac{1}{2}}
\]
is a norm continuous path of self-adjoint $J$-Fredholm operators.
The self-adjointness and boundedness follows for all $t\in[0,1]$ from
the spectral theorem applied to the function 
\[
x\mapsto x(1+x^2)^{-1/2}.
\]
So we need to prove the claims of continuity and $J$-Fredholmness.

For continuity, let $t,s\in [0,1]$, and apply 
\cite[Appendix A, Theorem 8]{CP1} to find
\[
\|F_{\diff_t}-F_{\diff_s}\|=\|\diff_t(1+\diff_t^2)^{-\frac{1}{2}}-
\diff_s(1+\diff_s^2)^{-\frac{1}{2}}\|\leq \|A_t-A_s\|
\]
proving continuity.

To prove the $J$-Fredholmness, let $t\in [0,1]$. Then 
\cite[Lemma $2.7$]{CP1} says that for $0<\varepsilon<1/4$ we have
$$F_{\D_t}-F_{\D_0}=B_\varepsilon(1+\D_0^2)^{-(1/2-\varepsilon)}$$
where $B_\varepsilon\in N$ and $\Vert B_\epsilon\Vert\leq
C(\epsilon)\Vert A_t-A_0\Vert$.
For $\varepsilon=1/4$ we get
\[
F_{\diff_t}-F_{\diff_0}=B_{1/4}(1+\diff^2_0)^{-1/4}.
\]
By \cite[Appendix, Lemma 6]{CP1}, and defining
$f(x)=1+\frac{x^2}{2}+\frac{x}{2}\sqrt{x^2+4}$, we have
$$(1+\diff_0^2)^{-1}\leq f(\Vert A_0\Vert)(1+\diff^2)^{-1}\in J,$$
and as $f(\Vert A_0\Vert)$ is scalar,
$$ (1+\D^2_0)^{-1/4}\leq  f(\Vert A_0\Vert)^{1/4}(1+\diff^2)^{-1/4}\in J.$$
Since $B_{1/4}$ is in $N$ we conclude
\[
\pi(F_{\diff_0})=\pi(F_{\diff_t}).
\]
At last for each $t\in[0,1]$,
\[
\pi(F_{\diff_t})\pi(F_{\diff_t})=\pi\big(\diff^2_t(1+\diff_t^2)^{-1}\big)=
\pi\big((1+\diff_t^2)(1+\diff_t^2)^{-1}\big)=\pi(1)
\]
so $\pi(F_{\diff_t})$ is invertible for all $t\in[0,1]$.

These considerations allow us to define spectral flow for such paths
of unbounded Fredholm operators. 
\begin{definition}[Unbounded Spectral Flow] Let $\{A_t\}_{t\in[0,1]}$
  be a norm continuous path of self-adjoint operators in $N$, and
  $(\A,H,\D)$ a von Neumann spectral triple relative to $(N,J)$.
The spectral flow of the  "continuous" path of unbounded selfadjoint
$J$-Fredholm operators $t\mapsto \diff+A_t$ is defined to be
\[
\textrm{\emph{sf}}\{\diff_t\}:=\textrm{\emph{sf}}\{F_{\diff_t}\}
\]
\end{definition}


\begin{prop}\label{unbound} Let $\{A_t\}_{t\in[0,1]}$
  be a norm continuous path of self-adjoint operators in $N$, and
  $(\A,H,\D)$ a von Neumann spectral triple relative to $(N,J)$. Let 
\[
p_1=\chi(F_{\diff+A_1})\conn
p_0=\chi(F_{\diff+A_0}).
\]
The spectral flow of the path $t\mapsto \diff+A_t$ only depends on
the end points $\diff+A_0$ and $\diff+A_1$ and is the class
\[
\begin{split}
\textrm{\emph{sf}}\{\diff_t\}=\textrm{\emph{sf}}\{F_{\diff_t}\}
&=[(1-p_1)\cap p_0]-[(1-p_0)\cap p_1]\\
&=\textrm{\emph{Ind}}_{(p_1\emph{-}p_0)}(p_1p_0)\in K_0(J).
\end{split}
\]
\end{prop}
\begin{proof}
Notice that
\[
\|\pi\big(\chi(F_{\diff_t})\big)-\pi\big(\chi(F_{\diff_s})\big)\|=
\|\chi\big(\pi(F_{\diff_t})\big)-\chi\big(\pi(F_{\diff_s})\big)\|=0
\]
for all $s,t\in[0,1]$ so by definition
\[
\textrm{sf}\{\diff_t\}=\textrm{sf}\{F_{\diff_t}\}=\big[(1-p_1)\cap
p_0\big]-\big[(1-p_0)\cap p_1\big]
\]
From this formula it is obvious that the spectral flow only depends
on the end points.
\end{proof}



\section{Kasparov Modules from Spectral Triples}



In this section we show that from any von Neumann spectral triple
$(\A,H,\D)$ relative to $(N,J)$, with $J$ $\sigma$-unital, we 
can construct a Kasparov module $(\mult_A,F_\D)\in \modu(A,J)$, where 
$\mult_A: A\to \ope(J)$ is left multiplication by elements in $A$. Defining 
$p_F=\frac{F_\D+1}{2}$ we get the class $[\mult_A,p_F]^1\in KK^1(A,J)$. We will
then show that for any unitary $u\in\A$, the unbounded spectral flow from $\D$ to
$u^*\D u$ is given by the Kasparov product
$[u]\grad_A[\mult_A,p_F]^1$. For an explanation of the terminology we refer to the 
appendix.

\subsection{Construction of a Kasparov Module}

Let $(\alg, H, \diff)$ be a von Neumann spectral triple relative to
$(N,J)$. Suppose that the norm-closed ideal $J$ is $\sigma$-unital. 

The ideal $J$ is a countably generated right Hilbert $J$-module when 
equipped with the inner product $\inner{x,y}=x^*y$ and
the action of $J$ from the right given by multiplication. Since $A$,
the norm closure of $\A$,
is represented in $N$, and $F_\D\in N$, we see that $F_\D\in \ope(J)$ 
and that there is a $*$-homomorphism 
$\mult_A : A\to \ope(J)$ given by left multiplication.

\begin{prop}\label{modulespec}
For all $a\in A$ the operators $[F_\diff,a]$, $a(1-F_\diff^2)$ and
$a(F_\diff-F_\diff^*)$ are in the norm closed ideal $J$. Thus, the pair 
$(\mult_A,F_\diff)$ is a Kasparov $A$-$J$-module.
\end{prop}
\begin{proof}
We have already noticed that $F_\diff=F_\diff^*$.
Let $a\in A$. Calculating modulo $J$
\[
aF_\diff^2=a\big(\diff^2(1+\diff^2)^{-1}\big)\sim
a\big(\diff^2(1+\diff^2)^{-1}+(1+\diff^2)^{-1}\big)=a
\]
So $a(F_\diff^2-1)\in J$ for all $a\in A$.

Let $a,b\in\alg$. We have
\[
[F_\diff,a]b=\diff\big[(1+\diff^2)^{-1/2},a\big]b+[\diff,a](1+\diff^2)^{-1/2}b
\]
As $[\diff,a]\in N$, see \cite[p.$456$]{CPRS2}, we have
\[
[\diff,a](1+\diff^2)^{-1/2}b\in J
\]
Thus we only need to show that
\[
\diff\big[(1+\diff^2)^{-1/2},a\big]b\in J
\]

Now, we employ integral formula, \cite[p.$8$]{Ped},
\[
(1+\diff^2)^{-1/2}=\frac{1}{\pi}\int_0^\infty
\lambda^{-1/2}(1+\diff^2+\lambda)^{-1}\,d\lambda.
\]
Denote the resolvent
$(1+\diff^2+\lambda)^{-1}$ by $R(\lambda)$. The provided
\[
\frac{1}{\pi}\int_{0}^\infty
\lambda^{-1/2}\diff\big[R(\lambda),a\big]b \,d\lambda
\]
is convergent in operator norm, it is equal to
\[
\diff\big[(1+\diff^2)^{-1/2},a\big]b.
\]
Applying some basic commutator identities yields
\[
\begin{split}
\diff\big[R(\lambda),a\big]b&=\diff
R(\lambda)[a,\diff^2]R(\lambda)b\\
&=\diff R(\lambda)[a,\diff]\diff R(\lambda)b+\diff
R(\lambda)\diff[a,\diff]R(\lambda)b.
\end{split}
\]
To establish the required norm estimates we require some inequalities.
The following inequalities can be proved using the spectral theorem for
unbounded operators, see \cite[Appendix A]{CP1},
\begin{enumerate}
\item $\|R(\lambda)\|=\|(1+\diff^2+\lambda)^{-1}\|\leq\frac{1}{1+\lambda}$
\item $\|\diff R(\lambda)\|=\|\diff(1+\diff^2+\lambda)^{-1}\|\leq
\frac{1}{2\sqrt{1+\lambda}}$
\item $\|\diff^2 R(\lambda)\|=\|\diff^2(1+\diff^2+\lambda)^{-1}\|\leq 1$
\end{enumerate}
for all $\lambda\in[0,\infty)$. Thus
\[
\frac{1}{\pi}\int_{0}^\infty
\lambda^{-1/2}\|\diff\big[R(\lambda),a\big]b\|
\,d\lambda\leq \frac{1}{\pi}\|b\|\|[a,\diff]\|
\int_{0}^\infty\lambda^{-1/2}\left(\frac{1}{4(1+\lambda)}+\frac{1}{1+\lambda}\right)
\,d\lambda<\infty.
\]
That is
\[
\diff\big[(1+\diff^2)^{-1/2},a\big]b=\frac{1}{\pi}\int_{0}^\infty
\lambda^{-1/2}\diff\big[R(\lambda),a\big]b \,d\lambda
\]
where the integral is convergent in operator norm. At last
\[
\begin{split}
\diff\big[R(\lambda),a\big]b&=\diff R(\lambda)ab -
[\diff,a]R(\lambda)b-a\diff R(\lambda)b\\
&=\diff R(\lambda)^{1/2}R(\lambda)^{1/2}ab- [\diff,a]R(\lambda)b
 -a\diff R(\lambda)^{1/2}R(\lambda)^{1/2}b\in J
\end{split}
\]
for all $\lambda\in[0,\infty)$ since all the terms are in $J$.
Thus we conclude that
\[
\diff\big[(1+\diff^2)^{-1/2},a\big]b\in J
\]
and thus that $[F_\diff,a]b\in J$ for all $a,b\in\alg$. By taking norm
limits $[F_\diff,a]b\in J$ for all $a,b\in A=\overline{\alg}$. 
\end{proof}

The argument used in the preceding proof is almost identical with the 
argument of S. Baaj and P. Julg used in \cite{Baaj} 
to build a bounded Kasparov module out of an unbounded one.



\subsection{The Pairing with $K_1(A)$ and Spectral Flow}



There is a certain case of unbounded spectral flow which is
particularly interesting. Suppose that $(\A,H,\D)$ is a unital 
 von Neumann spectral triple relative to $(N,J)$. 
Let $u\in\alg$ be unitary and consider the
path
\[
t\mapsto \diff_t:=(1-t)\diff+tu^*\diff u=\diff+ t[u^*,\diff]u
\]
The function $t\mapsto t[u^*,\diff]u$ is a continuous path of selfadjoint
elements in $N$, so we can calculate the spectral flow of the path $\diff_t$
via the transformation $\diff_t\mapsto F_{\diff_t}$.

\begin{lemma}\label{pinvert} Let $(\A,H,\D)$ be a unital 
 von Neumann spectral triple relative to $(N,J)$. 
Setting $p=\chi(F_\diff)$ and letting $u\in\A$ be unitary, we have 
$up-pu\in J$.
\end{lemma}
\begin{proof}
Polar decomposition of $F_\D$ gives
\[
F_\diff=(2p-1)|F_\diff|.
\]
So the image of $F_\D$ in the Calkin algebra is
\[
\pi(F_\diff)=\pi(2p-1)\pi(|F_\diff|)=\pi(2p-1)\pi(F_\diff^2)^{1/2}
=\pi(2p-1).
\]
It follows that
\[
2[u,p]-[u,F_\diff]=[u,(2p-1)-F_\diff]\in J.
\]
By Theorem \ref{modulespec} we have $[u,F_\diff]\in J$ so
$[u,p]\in J$ as claimed.
\end{proof}


\begin{prop}\label{uniflow}  Suppose that $(\A,H,\D)$ is a unital 
 von Neumann spectral triple relative to $(N,J)$, and $u\in\A$ is a unitary.
For the path $t\mapsto \diff_t$ from above we have
\[
\begin{split}
\emph{sf}\{\diff_t\}=\emph{sf}\{F_{\diff_t}\}
&=\partial\big[\pi(pup)+\pi(1-p)\big]\\
&=\textrm{\emph{Ind}}_{(p\emph{-}p)}(pup)
\end{split}
\]
where $p=\chi(F_\diff)$. From now on the spectral flow from $\diff$ 
to $u^*\diff u$ will be denoted by $\emph{sf}(\diff,u^*\diff u)$.
\end{prop}
\begin{proof}
From Theorem \ref{unbound}
\[
\textrm{sf}\{F_{\diff_t}\}=\big[(1-u^*pu)\cap
p\big]-\big[(1-p)\cap u^*pu\big]
\]
since
\[
\chi(F_{u^*\diff u})=\chi(u^*F_\diff u)=u^*\chi(F_\diff)u=u^*pu
\]
Now
\[
x\in \textrm{Ker}(u^*pu)\cap \textrm{Im}(p)\Leftrightarrow
\big(px=x\ \ \mbox{and}\ \  pux=0\big) 
\Leftrightarrow x\in \textrm{Ker}(pup)\cap
\textrm{Im}(p)
\]
and
\begin{align*}
x\in \textrm{Ker}(p)\cap \textrm{Im}(u^*pu)&\Leftrightarrow
\big(px=0\ \ \mbox{and}\ \  ux=pux\big)\\
&\Leftrightarrow \big(pu^*pux=0\ \ \mbox{and}\ \  ux=pux\big)
\Leftrightarrow ux\in
\textrm{Ker}(pu^*p)\cap \textrm{Im}(p).
\end{align*}
Thus
\[
(1-u^*pu)\cap p = N(pup)\cap p
\]
and
\[
u\big((1-p)\cap u^*pu\big)u^*=N(pu^*p)\cap p
\]
Since $N(pup+1-p)=N(pup)\cap p$ and $N(pu^*p+1-p)=N(pu^*p)\cap p$ we
conclude by Lemma \ref{connection} that
\[
\textrm{sf}\{F_{\diff_t}\}=\big[N(pup)\cap p\big]-\big[N(pu^*p)\cap
p\big]=\partial\big[\pi(pup)+\pi(1-p)\big]
\]
Remark that $\pi(pup)+\pi(1-p)$ is unitary in $N/J$ since $pu-up\in
J$.
\end{proof}

\begin{corollary}\label{idealexp}
Setting $p_F=\frac{F_\diff + 1}{2}$ we actually have
\[
\textrm{\emph{sf}}(\diff,u^*\diff u)=\partial[\pi(p_Fup_F+1-p_F)]
\]
\end{corollary}
\begin{proof}
In the proof of Lemma \ref{pinvert} we saw that $\pi(2p-1)=\pi(F_\diff)$ so
$\pi(p)=\pi(p_F)$ and the corollary follows easily.
\end{proof}

The last theorem of this section, which is the main theorem of the paper, expresses 
spectral flow from $\diff$ to $u^*\diff u$ in terms of a Kasparov product. 
We will need to use three different boundary maps namely 
\[
\begin{array}{ll}
\partial : K_1(N/J)\to K_0(J) & \\
\partial_{J\otimes\kom} : K_1\big(\calc(J\otimes\kom)\big)\to K_0(J\otimes\kom) & \textrm{and}\\
\partial_J : K_1\big(\calc(J)\big)\to K_0(J) & \\
\end{array}
\]
Note that for any $C^*$-algebra $B$, $\calc(B)$ denotes the Calkin
algebra 
$\ope(B)/B$.
Likewise we have the quotient maps 
\[
\begin{array}{llll}
\pi : N\to N/J &
\pi_{J\otimes\kom} : \ope(J\otimes\kom)\to \calc(J\otimes\kom)&\textrm{and} &
\pi_J : \ope(J)\to \calc(J)
\end{array}
\]

\begin{prop}\label{prodexp}
Suppose that the norm closed ideal $J$ is $\sigma$-unital and that the 
$C^*$-algebra $A=\overline{\alg}$ is separable.
Denote by $[\diff]=[\mult_A,p_F]^1$ the class in $KK^1(A,J)$ of 
the Kasparov module $(\mult_A,F_\diff)\in \modu(A,J)$ constructed in 
Theorem \ref{modulespec}.
Recall that $p_F=\frac{F_\diff+1}{2}$. 

Let $u\in\alg$ be unitary and
denote by $[u]$ its class in $K_1(A)$. Then we have the identity
\[
\textrm{\emph{sf}}(\diff,u^*\diff u)=\partial[\pi(p_Fup_F+1-p_F)]=
[u]\grad_A[\diff]
\]
\end{prop}
\begin{proof}
We start by stabilizing using the isomorphisms $K_1(A)\cong K_1(\Akom)$
and $KK^1(A,J)\cong KK^1(\Akom,J\otimes\kom)$.
In this way we obtain the classes
\[
[u\otimes e_{11}+e]\conn
[\mult_{\Akom}, p_F\otimes 1]^1
\]
in $K_1(\Akom)$ and $KK^1(\Akom,J\otimes\kom)$ respectively, where $e_{11}$
is a minimal projection in $\kom$ and $e=1-1\otimes e_{11}$. 
See \cite[Corollary $7.1.9$]{Wegge} and \cite[Corollary $17.8.8$]{Black}.
Thus, by Theorem \ref{pairing} the product is given by
\[
\begin{split}
[u]\grad_A[\diff]
&=\partial_{\Kom{J}}
\big[\pi_{\Kom{J}}\big(p_F\otimes 1 (u\otimes e_{11} + e) p_F\otimes 1+1-p_F\otimes 1\big)\big]\\
&=\partial_{\Kom{J}}\big[\pi_{\Kom{J}}
\big((p_Fup_F)\otimes e_{11} +p_F\otimes 1 -p_F\otimes e_{11}+1-p_F\otimes 1\big)\big]\\
&=\partial_{\Kom{J}}\big[\pi_{\Kom{J}}
\big((p_Fup_F+1-p_F)\otimes e_{11} + e\big)\big]
\end{split}
\]
in $K_0(\Kom{J})$. Recall that $\pi_J(p_F^2-p_F)=0$ since $[\mult_A,p_F]^1\in KK^1(A,J)$ 
and $A$ is unital.

But this is precisely the element
$\partial_J[\pi_J(p_Fup_F+1-p_F)]\in K_0(J)$ under
the isomorphism of $K_0(J)$ with $K_0(J\otimes\kom)$ 
\cite[Lemma $4.2.4$]{Higson}. Therefore the proof is finished if we can 
prove the identity 
\[
\partial_J[\pi_J(p_Fup_F+1-p_F)]=\partial[\pi(p_Fup_F+1-p_F)]
\]
To do so, let $x\in N$ be a norm-one lift of $\pi(p_Fup_F+1-p_F)\in N/J$. 
Then, as $N$ acts on $J$ by multiplication we have $N \subseteq \ope(J)$, so $x\in \ope(J)$ 
is too a norm-one lift 
of $\pi_J(p_Fup_F+1-p_F)\in \calc(J)$. Recalling the description of 
the boundary map using norm-one lifts given in \cite[Proposition $4.8.10$]{Higson}, 
the desired identity follows.  
\end{proof}



\section{$C^*$-Spectral Flow}


A problem with the construction of the Kasparov module in the last 
section is that it only works for $\sigma$-unital ideals $J$. For an 
arbitrary ideal in a von Neumann this may very well not be the case. 
Furthermore, if the ideal is $\sigma$-unital its $K$-theory is often 
simply $\rr$. When we can replace $J$ by a $\sigma$-unital 
$C^*$-algebra $B$, we not only ensure the existence of the $KK$-class, 
but can obtain stronger constraints on the values of the spectral flow.

\subsection{Basic Definitions}

Let $(\alg, H,\diff)$ be a unital von Neumann spectral triple relative
to $(N,J)$, and let  $A=\overline{\alg}$ be the norm closure 
of $\alg$. We assume that $A$ is a separable $C^*$-algebra.

Suppose that $B\subseteq J$ is a $\sigma$-unital $C^*$-algebra such that 
$(\mult_A,F_\diff)\in \modu(A,B)$, where $\mult_A : A\to \ope(B)$, thus 
in particular $A$ is supposed to act on $B$ by left-multiplication.

Note that $B$ is a countably generated right Hilbert $B$-module when equipped 
with the inner product $\inner{x,y}=x^*y$ for 
all $x,y\in B$ and the action of $B$ from the right given by multiplication.
The class  $[\mult_A, p_F]^1\in KK^1(A,B)$ is denoted by $[\diff_B]$ where 
$p_F=\frac{F_\D +1}{2}$.

Let $\partial_B : K_1\big(\calc(B)\big)\to K_0(B)$
be the boundary map, where $\calc(B)$ is the Calkin algebra $\ope(B)/B$. 
Let $\pi_B : \ope(B)\to \calc(B)$ denote the quotient map.

\begin{definition} Let $(\A,H,\D)$ and $B\subseteq J\subseteq N$ be as above.
We define the
$C^*$-\emph{spectral flow} as the quantity
\[
\textrm{\emph{sf}}_B(\diff,u^*\diff u)=\partial_B[\pi_B(p_Fup_F+1-p_F)]\in K_0(B)
\]
\end{definition}

The reason for supposing the existence of the Kasparov module class 
$[\diff_B]$ 
is that we want to describe the $C^*$-spectral flow using a Kasparov 
product. In fact we have

\begin{prop}\label{starproduct} Let $(\A,H,\D)$ be a von Neumann
  spectral triple relative to $(N,J)$.
Suppose that $B\subseteq J$ is a $\sigma$-unital $C^*$-algebra such 
that $(\mult_A,F_\diff)\in \modu(A,B)$ where $\mult_A : A\to \ope(B)$. 
Let $u\in \alg$ be unitary. The $C^*$-spectral flow from
$\diff$ to $u^*\diff u$ is equal to the product of 
$[\diff_B]\in KK^1(A,B)$ and the class of the
unitary $[u]\in K_1(A)$. That is
\[
[u]\grad_A[\diff_B]=\partial_B[\pi_B(p_Fup_F+1-p_F)]=
\textrm{\emph{sf}}_B(\diff,u^*\diff u)
\]
\end{prop}
\begin{proof}
The proof is similar to the one given in Theorem \ref{prodexp}.
\end{proof}

To justify the definition of $C^*$-spectral flow, we must show that
there exists a $\sigma$-unital
$C^*$-algebra $B$ such that $(\mult_A,F_\diff)$ is a Kasparov $A$-$B$-module. 

\begin{prop}  Let $(\A,H,\D)$ be a von Neumann
  spectral triple relative to $(N,J)$.
Let $\mathcal{B}$ be the smallest $C^*$-algebra in $\ope(H)$
containing the elements
\[
\begin{array}{cc}
F_\diff[F_\diff,a] & b[F_\diff,a]\\
F_\diff b[F_\diff,a] & a\varphi(\diff)
\end{array}
\]
for all $a,b\in \mathcal{A}$ and $\varphi\in C_0(\rr)$. Then
$\mathcal{B}$ 
is separable,
contained in $J$ and the pair 
$(\mult_A,F_\diff)$ is a Kasparov $A$-$\mathcal{B}$-module. 
In particular $\mathcal{B}$ is $\sigma$-unital. 
\end{prop}
\begin{proof}
Recall that $\alg$ is supposed to be unital. The $C^*$-algebra $C_0(\rr)$
is generated by the resolvent function $x\mapsto (i+x)^{-1}$
and the operator $(i+\diff)^{-1}$ is in $J$ so $a\varphi(\diff)$ is in
$J$ 
for all
$\varphi\in C_0(\rr)$. By Theorem \ref{modulespec}, $[F_\diff,a]\in J$
so 
all of the generators of
$\mathcal{B}$ are in $J$ and thus $\mathcal{B}\subseteq J$. Observe
that ${\mathcal B}$ is separable, and so $\sigma$-unital.

Now, clearly $A$ acts on $\mathcal{B}$ by multiplication from the left. 
Furthermore $F_\diff$ is in $\ope(\mathcal{B})$ since 
\[
\begin{array}{ccc}
1-F_\diff^2=(1+\diff^2)^{-1}\in \mathcal{B} & \textrm{and} & 
F_\diff\varphi(\diff)\in \mathcal{B}
\end{array}
\] 
for any $\varphi\in C_0(\rr)$. 

Proving that $(\mult_A,F_\diff)$ is a Kasparov 
$A$-$\mathcal{B}$-module is now straightforward.
\end{proof}



\subsection{The Relationship Between 
$C^*$-Spectral Flow and von Neumann 
Spectral Flow}



In this section we want to compare the $C^*$-spectral flow with the
von Neumann spectral flow. Let $(\alg,H,\diff)$ be a unital von Neumann 
spectral triple relative to $(N,J)$ and let
$u\in \alg$ be unitary. In Corollary \ref{idealexp} we found the expression
\[
\textrm{sf}(\diff,u^*\diff u)=\partial[\pi(p_Fup_F+1-p_F)]\in K_0(J)
\]
for the von Neumann spectral flow.

Let $B$ be a $\sigma$-unital $C^*$-algebra contained in 
$J$ such that $(\mult_A,F_\diff)$
is a Kasparov $A$-$B$-module. By definition we have the following expression
\[
\textrm{sf}_B(\diff,u^*\diff u)=\partial_B[\pi_B(p_Fup_F+1-p_F)]\in K_0(B)
\]
for the $C^*$-spectral flow.

These two notions should coincide in $K_0(J)$ when we apply the map
\[
i_*: K_0(B)\to K_0(J)
\]
induced by the inclusion $i : B\to J$.


\begin{lemma}\label{doublemult}  Let $(\alg,H,\diff)$ be a unital von Neumann 
spectral triple relative to $(N,J)$, and let $B$ be a 
$\sigma$-unital $C^*$-algebra contained in 
$J$ such that $(\mult_A,F_\diff)\in \modu(A,B)$ 
where $\mult_A : A\to \ope(B)$ is left-multiplication. 
The inclusion of $B$ in $\ope(H)$ can be extended to an injective 
$*$-homomorphism
\[
i : \ope(B)\to \ope(H)
\]
such that $i(T)(bx)=(Tb)x$ for all $T\in\ope(B)$, $b\in B$ and $x\in H$.
The image of the extension $i$ is contained in the double commutant of
$B\subseteq \ope(H)$. In particular $\ope(B)$ can be realized inside $N$.
\end{lemma}


\begin{proof}
Since $(\mult_A,F_\diff)$ is a Kasparov $A$-$B$-module and $A$ is unital,
we must have $1-F_\diff^2=(1+\diff^2)^{-1}\in B$. The image of
$(1+\diff^2)^{-1}\in \ope(H)$ is the domain of $\diff^2$ which is dense in $H$.
The representation of $B$ on $H$ by $i$ is thus seen to 
be non-degenerate.
Therefore, by \cite[Proposition $2.1$]{Lance}, the 
inclusion extends to $\ope(B)$
giving an injective $*$-homomorphism
\[
i : \ope(B)\to \ope(H)
\]
such that $i(T)(bx)=(Tb)x$ for all $T\in\ope(B)$, $b\in B$ and $x\in H$.

Let $S\in B'$ and let $T\in \ope(B)$. Suppose that 
$x\in H$ is of the form $x=by$
for some $b\in B$ and $y\in H$. Now
\[
i(T)Sby=i(T)bSy=(Tb)Sy=S(Tb)y=Si(T)by
\]
so $i(T)S=Si(T)$ on a dense subspace of $H$ and we conclude that 
$i(T)\in B''\subseteq N''=N$.
\end{proof}


\begin{prop}\label{coincide}  Let $(\alg,H,\diff)$ be a unital von Neumann 
spectral triple relative to $(N,J)$, and let $B$ be a 
$\sigma$-unital $C^*$-algebra contained in 
$J$ such that $(\mult_A,F_\diff)$
is a Kasparov $A$-$B$-module.
The von Neumann spectral flow coincides with the 
 $C^*$-spectral flow
under the homomorphism $i_* : K_0(B)\to K_0(J)$. More precisely for
$u\in\A$ unitary
\[
\textrm{\emph{sf}}(\diff,u^*\diff u)=
i_*\big(\textrm{\emph{sf}}_B(\diff,u^*\diff u)\big)
\]
\end{prop}
\begin{proof}
By Lemma \ref{doublemult} there are isometric maps 
\[
i: B\to J  \ \ \textrm{and}\ \   i: \ope(B)\to N 
\]
which allow us to also define the map (not necessarily injective)
\[
i: \calc(B)\to N/J
\]

Now, let $x\in \ope(B)$ be a norm-one lift of the unitary 
$\pi_B(p_Fup_F+1-p_F)\in \calc(B)$, then $i(x)\in N$ is a norm-one lift 
of the unitary $\pi(p_Fup_F+1-p_F)\in N/J$. By \cite[Proposition $4.8.10$]{Higson} 
we have 
\[
\begin{split}
& i_*\Big[\partial_B\big(\pi_B(p_Fup_F+1-p_F)\big)\Big]\\
&\qquad=i_*\left(\left[\begin{array}{cc}
xx^* & x(1-x^*x)^{1/2}\\
x^*(1-xx^*)^{1/2} & 1-x^*x
\end{array}\right]-
\left[\begin{array}{cc}
1 & 0 \\
0 & 0 
\end{array}\right]\right)\\
&\qquad=\left[\begin{array}{cc}
i(x)i(x)^* & i(x)\big(1-i(x)^*i(x)\big)^{1/2}\\
i(x)^*\big(1-i(x)i(x)^*\big)^{1/2} & 1-i(x)^*i(x)
\end{array}\right]-
\left[\begin{array}{cc}
1 & 0 \\
0 & 0 
\end{array}\right]\\
&\qquad=\partial[\pi(p_Fup_F+1-p_F)]
\end{split}
\]
Which is the desired identity. 
\end{proof}


\begin{corollary}\label{sepprod}  Let $(\alg,H,\diff)$ be a unital von Neumann 
spectral triple relative to $(N,J)$, and let $B$ be a 
$\sigma$-unital $C^*$-algebra contained in 
$J$ such that $(\mult_A,F_\diff)$
is a Kasparov $A$-$B$-module.
For $u\in\A$ unitary, the von Neumann spectral flow 
from $\diff$ to $u^*\diff u$ can be expressed in
terms of the Kasparov product of $[\diff_B]\in KK^1(A,B)$ and the class
$[u]\in K_1(A)$ of the unitary $u\in\alg$. 
More precisely
\[
\textrm{\emph{sf}}(\diff,u^*\diff u)=i_*([u]\grad_A[\diff_B])
\]
\end{corollary}
\begin{proof}
This follows immediately by Theorem \ref{coincide} and
Theorem \ref{starproduct}.
\end{proof}




\section{Numerical Spectral Flow}





Our aim in this section is to relate the von Neumann spectral flow to
the numerical spectral flow in semifinite von Neumann algebras 
studied in \cite{CPRS2,Ph,Ph1}.

In this section we let $N$ denote a semifinite von Neumann algebra
equipped with a fixed semifinite, faithful, normal trace $\tau$.
Furthermore, for any $*$-algebra $\fin\subseteq N$ we let $\fin^+$ denote
the $*$-algebra generated in $N$ by $\fin$ and the unit in $N$. When $\fin$ is
non-unital, we write the elements of $\fin^+$ as pairs $x+\lambda Id$, 
where $x\in\fin$ and $\lambda\in\cc$. 

\begin{definition}\label{taucompact}
Let $\mathcal{F}_N$ be the $*$-algebra in $N$ generated by the
projections 
$p$ with
finite trace, $\tau(p)<\infty$. By \cite[Section 1.8]{Fack}, 
$\mathcal{F}_N$
is an ideal in $N$.
The \emph{$\tau$-compact operators}, $\kom_N$, is the norm-closure of 
$\mathcal{F}_N$.
\end{definition}

Let $(\A,H,\D)$ be a semifinite spectral triple relative to $(N,\tau)$ as defined 
in \cite[Definition $2.1$]{CPRS2}. Notice that $(\A,H,\D)$ is a von Neumann spectral 
triple relative to $(N,\kom_N)$ in an obvious way. For semifinite spectral triples, 
spectral flow is defined as a real number, whereas our methods produce a class
in $K_0(\kom_N)$. The problem is solved by 
establishing a homomorphism $\tau_* : K_0(\kom_N)\to \rr$
by means of the  trace $\tau$. The existence and nature of such a 
homomorphism is of course well known,
but as the link to the semifinite case is very important we will 
carry out the details.


\begin{lemma}
Let $n\in\nn$. For each finite set of elements 
$\{x_1,\ldots,x_m\}\subseteq M_n(\fin_N)$
there is a projection $p\in M_n(\fin_N)$ with $px_i=x_i$ for all 
$i\in\{1,\ldots,m\}$.
The projection $p$ is called a \emph{local unit} for $\{x_1,\ldots,x_m\}$.
\end{lemma}
\begin{proof} 
For any finite set of projections
$\{p_1,\ldots,p_m\}$ in $\fin_N$,
the inequality $\sup\{p_1,\ldots,p_m\}\leq p_1+\ldots+p_m$ holds so 
$\sup\{p_1,\ldots,p_m\}\in \fin_N$.
Furthermore, for each $i\in \{1,\ldots,m\}$ we have 
$p_i\leq \sup\{p_1,\ldots,p_m\}$ yielding
$\sup\{p_1,\ldots,p_m\}p_i=p_i$, so $\sup\{p_1,\ldots,p_m\}$ is a local unit for 
$\{p_1,\ldots,p_m\}$. To obtain the desired property for $\fin_N$, 
note that each element in $\fin_N$ is a complex polynomial of finite degree, where 
the variables are projections with finite trace.

Now, let $n\in\nn$ and fix a finite set of matrices
$\{x_1,\ldots,x_m\}\subseteq M_n(\fin_N)$. Choose a projection $p\in \fin_N$,
such that $px_i^{kl}=x_i^{kl}$ for all $i\in \{1,\ldots,m\}$ and
$k,l\in \{1,\ldots n\}$, where $x_i^{kl}$ is the matrix entry corresponding 
to row $k$ and column $l$. Then obviously $\diag(p,\ldots p)x_i=x_i$ for all
$i\in \{1,\ldots, m\}$ as desired.
\end{proof}


\begin{lemma}\label{isomfin}
For each $n\in \nn$ the $*$-algebra $\fin_N$ is stable under the 
holomorphic functional calculus.
That is, for $x\in M_n(\fin_N)$ and for 
$f$ a holomorphic function in a neighborhood of the spectrum of 
$x$ in $M_n(\kom_N)$ with $f(0)=0$ 
we have $f(x)\in M_n(\fin_N)$. In other words, $\fin_N$ equipped with the $C^*$-norm 
from $\kom_N$ becomes a pre-$C^*$-algebra. In particular it has a well-defined 
$K$-theory and, by \cite[Proposition $3$]{Connes}, 
the inclusion $i: \fin_N\to \kom_N$ induces an isomorphism 
$i_* : K_0(\fin_N)\to K_0(\kom_N)$.
\end{lemma}
\begin{proof} We employ the technique of \cite[Proposition $4$]{ar}.
First, notice that $f(x)\in M_n(\kom_N)$ because $M_n(\kom_N)$ is a 
$C^*$-algebra. Now,
for a closed curve $\gamma$ winding once around the spectrum of $x$ in $M_n(\kom_N)$ 
not touching $0$, the identity
\[
f(x)=\frac{1}{2\pi i}\int_{\gamma}f(\lambda)(\lambda-x)^{-1}d\lambda
\]
is valid. Let $p$ be a local unit for $x$, let $\lambda$ be in the resolvent of $x$
and check that
\[
(1-p)=(1-p)(x-\lambda)(x-\lambda)^{-1}=-\lambda(1-p)(x-\lambda)^{-1}.
\]
Thus for $\lambda\neq 0$ we have
\[
(1-p)(x-\lambda)^{-1}=-\frac{1}{\lambda}(1-p)
\]
This enables us to calculate
\[
\begin{split}
(1-p)f(x)&=\frac{1}{2\pi i}\int_{\gamma}f(\lambda)(1-p)(\lambda-x)^{-1}d\lambda\\
&=\frac{1}{2\pi i}\int_{\gamma}\frac{-f(\lambda)}{\lambda}(1-p)d\lambda \\
&=(1-p)f(0)=0
\end{split}
\]
It follows that $pf(x)=f(x)$. As $M_n(\fin_N)$ is an ideal in 
$M_n(\kom_N)$ and $p\in M_n(\fin_N)$ we conclude that $f(x)\in M_n(\fin_N)$ as desired.
\end{proof}


\begin{prop}\label{homokom}
There is a homomorphism $\tau_* : K_0(\kom_N)\to \rr$ given by
\[
\tau_*\big([x+\lambda Id]-[y+\mu Id]\big)=\tau_n(x)-\tau_n(y)
\]
for each pair of projections $x+\lambda Id,\ y+\mu Id\in M_n(\fin_N^+)$ with
$[\lambda]=[\mu]$ in $K_0(\cc)$. Remark that $\tau_n=\tau\otimes \emph{Tr}$ 
on the algebraic tensor product $\fin_N\otimes M_n(\cc)=M_n(\fin_N)$ where 
$\emph{Tr}$ is the canonical trace on $M_n(\cc)$.
\end{prop}
\begin{proof}
Define $\hat\tau : \fin_N^+\to \rr$ by $\hat\tau(x+\lambda
Id)=\tau(x)$ then 
$\hat\tau$
satisfies the relation $\hat\tau(u^*xu)=\hat\tau(x)$ for all unitaries
$u\in \fin_N^+$. 
Indeed, write $u=v+\alpha Id$, where $v\in \fin_N$ and $\alpha\in \cc$, then
\[
\overline{\alpha}\alpha=1\conn
v^*v+v^*\alpha+v\overline{\alpha}=0=vv^*+v^*\alpha+v\overline{\alpha}
\]
Now, we simply calculate
\[
\begin{split}
(v^*+\overline{\alpha}Id)(x+\lambda Id)(v+\alpha Id)
&=(v^*xv+v^*x\alpha+v^*\lambda v+ v^*\lambda\alpha+\overline{\alpha}xv+
x+\overline{\alpha}\lambda v+\lambda Id)\\
&=(v^*xv+v^*x\alpha+\overline{\alpha}xv + x+\lambda Id)
\end{split}
\]
thus applying our extended $\hat\tau$ yields
\[
\hat\tau\big((v^*+\overline{\alpha}Id)(x+\lambda Id)(v+\alpha Id)\big)=
\tau(v^*xv+v^*x\alpha+\overline{\alpha}xv + x)=
\tau(x)
\]
Now, clearly, there is a well-defined homomorphism 
$\tau_* : K_0(\fin_N^+)\to \rr$ given by
$\hat\tau([x+\lambda Id]-[y+\mu Id])=\tau_n(x)-\tau_n(y)$ 
for each pair of projections
$(x+\lambda Id),(y+\mu Id)\in M_n(\fin_N^+)$. 
Since $K_0(\fin_N)$ is the kernel of the homomorphism
$\pi_* : K_0(\fin_N^+)\to K_0(\cc)$ induced by the projection 
$\pi : \fin_N^+\to \cc$ we get the
desired map by restriction and a reference to Lemma \ref{isomfin}.
\end{proof}


\begin{prop}\label{fintrace}
Let $p$ be a projection in $M_n(\kom_N)$, then actually $p\in M_n(\fin_N)$.
\end{prop}
\begin{proof}
Since $M_n(\fin_N)$ is dense in $M_n(\kom_N)$, there is a positive element
$e\in M_n(\fin_N)$ such that $\|e-p\|<\frac{1}{24}$. In particular
$\|e\|<2$. The estimate
\[
\|e^2-e\|\leq\|e(e-p)\|+\|(e-p)p\|+\|p-e\|<\frac{1}{4} 
\]
shows that $e$ is almost a projection. It follows that $1/2\notin \spec(e)$, 
creating a gap in the spectrum of $e$.
There is thus an $\varepsilon >0$
such that
\[
\spec(e)\subseteq [0,1/2-\varepsilon]\cup [1/2+\varepsilon,5/4]
\]
and the function $f:\rr/\{\frac{1}{2}\}\to\rr$ given by
\[
f(t)=\left\{\begin{array}{cc}
0 & \qquad t<\frac{1}{2} \\
1 & \qquad t>\frac{1}{2}
\end{array}\right.
\]
is holomorphic on a neighborhood of $\spec(e)$ with $f(0)=0$. 
By Lemma \ref{isomfin} the projection $f(e)$ is in $M_n(\fin_N)$. Moreover,
\[
\sup\{|f(t)-t|\,|\, t\in \spec(e)\}\leq \sup\{1/2-\varepsilon,1/4\}
\]
so $\|f(e)-e\| <\frac{1}{2}$. This gives us the inequality
\[
\|p-f(e)\|\leq \|p-e\|+\|e-f(e)\|<1
\]
but then $p$ and $f(e)$ must be equivalent, i.e.
there exist a unitary $u$ in $M_n(\kom_N^+)$ such that $u^*f(e)u=p$ by
\cite[Proposition $4.1.7$]{Higson}.
The proof is finished by recalling that $M_n(\fin_N)$ is an ideal in $M_n(N)$.
\end{proof}



\begin{definition}\cite{Ph,Ph1}
Let $\{B_t\}$ be a path of selfadjoint operators in $N$ such that
$\pi(B_t)\in N/\kom_N$ is invertible for all $t$. Let
$0=t_0<t_1<\ldots<t_n=1$ be a partition of $[0,1]$ such that for each
$i\in\{1\ldots,n\}$ we have
\[
\|\pi\big(\chi(B_t)\big)-\pi\big(\chi(B_s)\big)\|<1/2
\]
for all $t,s\in [t_{i-1},t_i]$. Recalling that all projections in 
$\kom_N$ have 
finite trace by Theorem \ref{fintrace}, we define the 
semifinite spectral flow of the path $\{B_t\}$ as the real number
\[
\textrm{\emph{sf}}\{B_t\}=
\sum_{i=1}^n\Big(\tau\big(N(p_i)\cap p_{i-1}\big)-
\tau\big(N(p_{i-1})\cap p_i\big)\Big)
\]
where $p_i=\chi(B_{t_i})$. 
\end{definition}

\begin{prop}
The semifinite spectral flow of the path $\{B_t\}$ from above can be expressed as
\[
\emph{sf}\{B_t\}=
\tau\big(N(p_n\ldots p_0)\cap p_0\big)-
\tau\big(N(p_0\ldots p_n)\cap p_n\big)
\]
Moreover it is independent of the partition chosen and is invariant under homotopies 
of the path $\{B_t\}$ keeping the endpoints fixed.
\end{prop}
\begin{proof}
This follows immediately by applying our homomorphism
$\tau_* : K_0(\kom_N)\to \rr$ from Theorem \ref{homokom} to the results in
Corollary \ref{beginning}, Lemma \ref{conclusion} and
Lemma \ref{homotopy}
recalling that each projection $p\in\kom_N$ has finite trace by Theorem \ref{fintrace}
\end{proof}

\begin{definition}\cite{CP1,CP2,Ph,Ph1} Let $(\alg,H,\diff)$ be a 
unital semifinite spectral triple relative
to $(N,\tau)$. Suppose that the norm closure $A=\overline{\alg}$ 
of $\alg$, is a separable $C^*$-algebra.
For each path of selfadjoint operators $\{A_t\}$ in $N$, we
define the  semifinite spectral flow of the 
path $t\mapsto \diff + A_t:=\diff_t$
as the real number $\textrm{sf}\{\diff_t\}:=\textrm{sf}\{F_{\diff_t}\}$.
\end{definition}

We can now state our main theorem relating the three different notions
of spectral flow we have discussed.

\begin{prop}\label{maintheorem} 
Let $(\alg,H,\diff)$ be a unital semifinite spectral 
triple relative
to $(N,\tau)$. Suppose that the norm closure $A=\overline{\alg}$ 
of $\alg$, is a separable $C^*$-algebra.
Let $u\in\alg$ be unitary. Set 
$\diff_t = (1-t)\diff + tu^*\diff u=\diff +tu^*[\diff,u]$,
then the unbounded semifinite spectral flow of the path $t\mapsto\diff_t$ is
given by
\[
\textrm{\emph{sf}}\{\diff_t\}=
\tau\big(\partial[\pi(pup+1-p)]\big)=
\tau\big(N(pup+1-p)\big)-\tau\big(N(pu^*p+1-p)\big)
\]
where $\tau_* : K_0(\kom_N)\to \rr$ is the homomorphism from Theorem \ref{homokom} and
$p=\chi(F_\diff)$.
In addition there exists a separable $C^*$-algebra $\mathcal{B}\subseteq \kom_N$ and a 
class $[\diff_{\mathcal{B}}]\in KK^1(A,\mathcal{B})$ such that 
\[
\textrm{\emph{sf}}\{\diff_t\}=\tau_*\big(i_*([u]\grad_A[\diff_B])\big)
\]
where $i :\mathcal{B} \to \kom_N$ is the inclusion and 
$[u]\in K_1(A)$ is the class of the unitary.
\end{prop}
\begin{proof}
This follows immediately by applying our $\tau_* :K_0(\kom_N)\to \rr$ from
Theorem \ref{homokom} to both sides of the equalities in Theorem \ref{uniflow}
and in Corollary \ref{sepprod}, keeping in mind that each projection
in $\kom_N$ has finite trace by Theorem \ref{fintrace}.
\end{proof}

Theorem \ref{maintheorem} shows that semifinite spectral triples
represent $KK$-classes in a precise sense. While this is really proved
here only for odd spectral triples, the discussion in \cite{prs} and
some simple adaptations of these proofs show that such a
representation theorem is also true in the even case.



\section{Appendix on Kasparov Products}




In this appendix we give explicit forms for odd pairings in
$KK$-theory. In order to do this, we need to recall some basic
definitions and results.

\begin{definition}
Let $A$ and $B$ be $\zz_2$-graded $C^*$-algebras.
A Kasparov $A$-$B$-module is a pair $(\psi,V)$ consisting of a 
graded $*$-homomorphism $\psi : A\to \ope(E)$, with $E$ a countably generated, graded  
right Hilbert $B$-module, together with an odd operator $V\in \ope(E)$ such that 
\begin{enumerate}
\item $\psi(a)(V^2-1)\in \kom(E)$
\item $\psi(a)(V-V^*)\in \kom(E)$
\item $[V,\psi(a)]\in \kom(E)$
\end{enumerate}
for all $a\in A$. 
The set of Kasparov $A$-$B$-module is denoted by $\modu(A,B)$. An element 
$(\psi,V)\in \modu(A,B)$ is called degenerate when 
$a(V^2-1)=a(V-V^*)=[V,a]=0$.
\end{definition}

The set $\modu(A,B)$ becomes the even $KK$-theory group $KK(A,B)$ 
when equipped with direct sum 
and the equivalence relation $\sim_{oh}$ generated by operator
homotopy, unitary equivalence and addition of degenerate elements.
Unitary equivalence is denoted by $\sim_u$. The class represented by 
the pair $(\psi,V)\in \modu(A,B)$ is denoted by $[\psi,V]\in KK(A,B)$. 
\cite[Proposition $17.3.3$]{Black}.

To define the odd $KK$-theory group we introduce the Clifford algebra 
$\cc_1$, that is the $C^*$-algebra $\cc\oplus\cc$ 
equipped with the standard odd grading, and we set $KK^1(A,B)=KK(A,B\grad\cc_1)$ 
where $\grad$ denotes the graded tensor product as defined in \cite[Chapter $14.4$]{Black}. 

For ungraded $C^*$-algebras $A$ and $B$ there is a description of 
the odd $KK$-theory using extensions of $C^*$-algebras. More precisely 

\begin{prop}\cite[Proposition $17.6.5$]{Black}\label{extisom}
There is an isomorphism 
\[
\textrm{\emph{Ext}}^{-1}(A,B\otimes\kom)\cong KK^1(A,B\otimes\kom)
\]
An invertible extension given by the $*$-homomorphism
$\psi : A\to \ope(\Bkom)$ and the element
$p\in \ope(\Bkom)$ that is with Busby-invariant
$\tau : a\mapsto \pi\big(p\psi(a)p\big)\in \calc(\Bkom)$
is mapped to the Kasparov $A$-$(\Bkom)$-module
$\big(\psi\grad 1, (2p-1)\gradodd\big)\in \modu\big(A,(\Bkom)\grad\cc_1\big)$ 
with $\psi\grad 1 : A\to \ope\big((\Bkom)\grad\cc_1\big)$ and 
$\varepsilon=(1,-1)\in \cc_1$
\end{prop}

With this in mind we will employ the notation
$[\psi,p]^1$ 
for the 
class $[\psi\grad 1,(2p-1)\gradodd]\in KK^1(A,\Bkom)=
KK\big(A,(\Bkom)\grad\cc_1\big)$ 
where $\psi : A\to \ope(\Bkom)$ and 
$p\in \ope(\Bkom)$ have the properties 
\begin{enumerate}
\item $\psi(a)(p^2-p)\in \Bkom$
\item $\psi(a)(p-p^*)\in \Bkom$
\item $[p,\psi(a)]\in \Bkom$
\end{enumerate}
for all $a\in A$.

Let $A$, $B$ and $D$ be graded $C^*$-algebras. Suppose that $A$ and
$D$ 
are separable 
and that $B$ is $\sigma$-unital. 
A fundamental property of $KK$-theory is the existence of a bilinear 
associative product 
\[
\grad_A : KK^i(D,A)\times KK^j(A,B)\to KK^{i+j}(D,B)
\]
The aim of this appendix is to give a concrete description of a certain 
instance of this product namely the one between $K_1(A)=KK^1(\cc,A)$ 
and $KK^1(A,B)$, \cite[Chapter $18$]{Black}.








We will need to quote a couple of results. First of all, since the aim is 
to form products with $K$-theory we will use the isomorphism of $K$-theory 
with $KK$-theory. 


\begin{lemma}\label{isomk}
Let $A$ be an ungraded $C^*$-algebra. 
The groups $KK^1(\cc,\Akom)$ and $K_1(\Akom)$ are isomorphic.
The isomorphism is given by
\[
[\mult_\cc,p]^1\mapsto \partial[\pi(p)]
\]
where $\mult_\cc : \cc\to \ope(\Akom)$ is left multiplication by the 
complex numbers, $\pi : \ope(\Akom)\to \calc(\Akom)$ is the quotient map and 
$\partial : K_0\big(\calc(\Akom)\big)\to K_1(\Akom)$ is the boundary map,
\cite[Proposition $17.5.7$]{Black}.

>From the definition of $\partial$ it follows that, for each class
$[u]\in K_1(\Akom)$, there exists a selfadjoint 
$q\in \ope(\Akom)$ with $\|q\|\leq 1$
such that $[u]=[\exp(2\pi i q)]$, \cite[Proposition $12.2.2$]{Roerdam}.


Likewise the groups $KK(\cc,\Akom)$ and  $K_0(\Akom)$ are isomorphic.
The isomorphism is given by
\[
[\mult_\cc, V]\mapsto \partial[\pi(T)]
\]
where $\mult_\cc : \cc \to \ope\big((\Akom)\oplus(\Akom)\big)$, the element 
$V\in \ope\big((\Akom)\oplus(\Akom)\big)$ is the matrix 
\[
V=\left(
\begin{array}{cc}
0 & T^*\\
T & 0
\end{array}
\right)
\]
and $\partial : K_1\big(\calc(\Akom)\big)\to K_0(\Akom)$ is the boundary map.
Note that the grading on $(\Akom)\oplus(\Akom)$ is given by $\gamma = \left(
\begin{array}{cc}
0 & 1\\
1 & 0
\end{array}
\right)$, \cite[Proposition $17.5.5$]{Black}.
\end{lemma}


Let $A$, $B$ and $D$ be graded $C^*$-algebras, with $A$ and $D$ separable and $B$ 
$\sigma$-unital.
The Kasparov product can be constructed using the notion of connections. 
Let $(\psi_1,V_1)\in \modu(D,A)$ with $\psi_1: D\to \ope(E_1)$ 
and let $(\psi_2,V_2)\in \modu(A,B)$ with $\psi_2: A\to \ope(E_2)$.
We can form the graded interior tensor product 
$E=E_1\grad_{\psi_2} E_2$ in the sense of \cite[Chapter $14.4$]{Black}. For each $x\in E_1$ there 
is a map $T_x\in \ope(E_2,E)$ such that $T_x(y)=x\grad y$ for all 
$y\in E_2$, \cite[Proposition $4.6$]{Lance}.


\begin{definition}\label{appconn}
An odd operator $F\in \ope(E)$ is called a $V_2$-connection for $E_1$ if, 
for any homogeneous $x\in E_1$, we have
\[
T_x V_2 - (-1)^{\partial x}FT_x\in \kom(E_2,E)
\]
where $\partial x$ denotes the degree of $x$ in $E_1$.
\end{definition}


Now we are ready to state the most important background result. It gives a 
concrete description of the product under an assumption on commutators. 
Later on the $C^*$-algebra $D$ is going to be the complex numbers so the 
assumption will be trivially satisfied.


\begin{prop}\cite[Proposition $18.10.1$]{Black}\label{evenprod}
Let $x=(\psi_1,V_1)\in \modu(D,A)$ with $V_1=V_1^*$ and $\|V_1\|\leq 1$.
Let $y=(\psi_2,V_2)\in \modu(A,B)$. Let $F$ be a $V_2$-connection for $E_1$.
Set $E=E_1\grad_{\psi_2}E_2$, $\psi=\psi_1\grad 1 : A\to \ope(E)$ and
\[
V=V_1\grad 1 + \big((1-V_1^2)^{1/2}\grad 1\big) F
\]
If $[V_1\grad 1,\psi(a)]\in \kom(E)$ for all 
$a\in A$, then $z=(\psi,V)\in \modu(D,B)$ is
operator homotopic to the Kasparov product of $x$ and 
$y$, i.e. $[x]\grad_A[y]=[z]$
in $KK(D,B)$.
\end{prop}


To form the product in $KK^1$ we need to be able to 
move the Clifford algebra from 
the second coordinate to the first. This is accomplished
 by the following lemma.


\begin{lemma}\label{movecliff}
Let $A$ and $B$ be ungraded $C^*$-algebras. There is a group isomorphism
\[
\varphi : KK^1(A,\Bkom)=KK\big(A,(\Bkom)\grad\cc_1\big)\to
KK(A\grad\cc_1,\Bkom)
\]
such that 
\[
\varphi[\psi,p]^1=\left[\sigma, \left(\begin{matrix}0&2p-1\\
 2p-1&0\end{matrix}\right)\right]
\]
where the graded $*$-homomorphism $\sigma : A\grad\cc_1\to 
\ope\big((\Bkom)\oplus (\Bkom)\big)$ is 
given by
\[
\sigma(a,-a)=\left(\begin{matrix}0&-i\psi(a)\\
i\psi(a)&0\end{matrix}\right)
\]
and
\[
\sigma(a,a)=\left(\begin{matrix}\psi(a)& 0\\
0&\psi(a)\end{matrix}\right)
\]
The grading on $(\Bkom)\oplus (\Bkom)$ is given by the grading operator 
$\gamma = \left(\begin{matrix}0& 1\\
1&0\end{matrix}\right)$, \cite{K}.
\end{lemma}




\subsection{Product between $K_1$ and $KK^1$}




The starting point is a translation of Theorem \ref{evenprod} suited to handle the 
odd case.

\begin{prop}\label{prod}
Let $A$, $B$ and $D$ be ungraded $C^*$-algebras, with $A$ and $D$ separable and
$B$ $\sigma$-unital. Let $[x]$ be a class in $KK^1(D,\Akom)$. By 
Theorem \ref{extisom} we can assume that $[x]$ is represented by the 
Kasparov module 
\[
x=(\psi_1\grad 1, (2q-1)\gradodd)\in\modu\big(D,(\Akom)\grad\cc_1\big)
\]
where $\psi_1\grad 1 : D\to \ope(E_1)$, with $E_1=(\Akom)\grad\cc_1$. 
By \cite[Proposition $17.4.3$]{Black} we may assume that $q=q^*$ 
and that $\|q\|\leq 1$.

Let $[y]$ be a class in $KK\big((\Akom)\grad\cc_1,B\big)\cong KK^1(\Akom,B)$. 
See Lemma \ref{movecliff}. Suppose that $[y]$ is represented by the module 
\[
y=(\psi_2,V_2)\in \modu\big((\Akom)\grad\cc_1, B\big)
\]
with $\psi_2 :  (\Akom)\grad\cc_1\to \ope(E_2)$.

Set $E=E_1\grad_{\psi_2}E_2$ and
$\psi=(\psi_1\grad 1)\grad 1$. Let $F\in \ope(E)$ be a
$V_2$-connection for $E_1$. Define
\[
V=-\big(\cos(\pi q)\gradodd\big)\grad 1 + \Big(\big(\sin(\pi q)\grad
1\big)\grad 1\Big)F\in \ope(E)
\]
Suppose that
\[
\Big[\big(\cos(\pi q)\grad\varepsilon\big)\grad 1,\psi(d)\Big]\in
\kom(E)
\]
for all $d\in D$. Then $(\psi,V)$ is a Kasparov $D$-$B$-module
which is operator homotopic to the Kasparov product
of $x$ and $y$. That is $[\psi,V]=[x]\grad_{\Akom}[y]$.
\end{prop}
\begin{proof}
Let $d\in D$. Remark that $\psi_1(d)(q^2-q)\in \Akom$, thus modulo
$\Akom$ we have
\[
\begin{split}
\psi_1(d)\cos(\pi q)
&=\psi_1(d)\sum_{k=1}^\infty (-1)^k\frac{(\pi q)^{2k}}{(2k)!}+\psi_1(d)\\
&\sim\psi_1(d)\sum_{k=1}^\infty (-1)^k\frac{\pi^{2k}}{(2k)!}q+\psi_1(d)\\
&=\psi_1(d)\big(\cos(\pi)q-q+1\big)\\
&=-\psi_1(d)(2q-1)
\end{split}
\]
for all $d\in D$. It follows that $x$ is a compact pertubation of
\[
(\psi_1\grad 1, -\cos(\pi q)\gradodd\big)\in 
\modu\big(D,(\Akom)\grad\cc_1\big)
\]
so they determine the same class in $KK\big(D,(\Akom)\grad\cc_1\big)$.

By assumption the last module fulfils the conditions of Theorem \ref{evenprod}
so the theorem is proved if
\[
\Big(1-\big(\cos(\pi q)\gradodd\big)^2\Big)^{1/2}=\sin(\pi q)\grad 1
\]
but this is clear since we supposed that $\|q\|\leq 1$ and 
$q=q^*$ so $\spec(q)\subseteq [-1,1]$, 
a fact which yields the positivity of $\sin(\pi q)\grad 1$.
\end{proof}



\begin{prop}\label{pairing}
Suppose that $A$ and $B$ are ungraded $C^*$-algebras, with $A$ separable 
and $B$ $\sigma$-unital. 
Let $[u]\in K_1(\Akom)$. The isomorphism from Lemma \ref{isomk} sends $[u]$ 
to a class $[\mult_\cc, q]^1\in KK^1(\cc,\Akom)$
where $q=q^*$ and $\|q\|\leq 1$. In particular $[u]=[\exp(2\pi i q)]$ so
without loss of generality we can assume that $u=\exp(2\pi i q)$.
Let $y=[\psi,p]^1\in KK^1(\Akom,\Bkom)$ and assume that 
$\psi : \Akom\to\ope(\Bkom)$ is non-degenerate. Then the product
$[u]\grad_{\Akom}y$ is equal to
\[
\partial\Big[\pi\big(p\psi(u)p+(1-p)\big)\Big]\in K_0(B\otimes\kom)
\]
where $\psi$ is extended to $(\Akom)^+$ and 
$\partial : K_1\big(\calc(\Bkom)\big)\to K_0(\Bkom)$ is the boundary map.
\end{prop}

\begin{proof}
Applying the isomorphism
\[
\varphi : KK^1(\Akom,\Bkom)\to KK\big((\Akom)\grad\cc_1,\Bkom\big)
\]
from Lemma \ref{movecliff} to $y$ we get
\[
\begin{split}
\varphi y=\varphi[\psi, p]^1=[\sigma,V_2]
\end{split}
\]
with $V_2= \left(\begin{matrix}0& 2p-1\\
 2p-1&0\end{matrix}\right)$ and 
$\sigma : (\Akom)\grad\cc_1\to \ope\big((\Bkom)\oplus(\Bkom)\big)$
given by 
\[
\begin{array}{ccc}
\sigma(a,-a)=\left(\begin{matrix}0&-i\psi(a)\\
i\psi(a)&0\end{matrix}\right)&
\textrm{and} &
\sigma(a,a)=\left(\begin{matrix}\psi(a)& 0\\
0&\psi(a)\end{matrix}\right)
\end{array}
\]
which thus canonically represents $y$ in
$KK\big((\Akom)\grad\cc_1,\Bkom\big)$.

Recall that $[\mult_\cc, q]^1\in KK^1(\cc, \Akom)$ is notation for the class 
$[x]\in KK\big(\cc, (\Akom)\grad\cc_1\big)$ represented by the module 
\[
x=\big(\mult_\cc, (2q-1)\gradodd\big)\in \modu\big(\cc, (\Akom)\grad\cc_1\big)
\]

We are now in position to form the product
$z=[u]\grad_{\Akom}y=[x]\grad_{\Akom}\varphi y$.
Set $E_1 = (\Akom)\grad\cc_1$, $E_2 = (\Bkom)\oplus(\Bkom)$ and
$E=E_1\grad_\sigma E_2$. Recall that the grading on $(\Bkom)\oplus(\Bkom)=E_2$ is given by the
grading operator $\gamma = \left(\begin{array}{cc}0&1\\
1&0\end{array}\right)$.
Since $\psi$ is assumed to be 
non-degenerate, $\sigma$ is also non-degenerate and there is an even unitary
isomorphism
\[
\begin{split}
w : E_1\grad_\sigma E_2&=\big((\Akom)\grad\cc_1\big)
\grad_\sigma\big((\Bkom)\oplus(\Bkom)\big)\\
&\to(\Bkom)\oplus(\Bkom)=E_2
\end{split}
\]
given by
\[
w(x_1\grad_\sigma x_2)=\sigma(x_1)x_2\qquad x_1\in E_1\,,\,x_2\in
E_2
\]

Let $x\in E_1$ be homogeneous. Clearly $wT_x = \sigma(x)$ so
\[
wT_xV_2-(-1)^{\partial x}V_2wT_x=\big[\sigma(x),V_2\big]\in
\kom(E_2)
\]
Thus $w^*V_2w\in E$ is a $V_2$-connection for $E_1$. By Theorem \ref{prod} we
can represent the product $z$ by the module 
$(M_\cc,V)$ where
\[
V=-\big(\cos(\pi q)\gradodd\big)\grad 1 + \Big(\big(\sin(\pi q)\grad
1\big)\grad 1\Big)w^*V_2w\in \ope(E)
\]
But $(\mult_\cc,V)\sim_u (M_\cc,wVw^*)$ so
actually
\[
z=[\mult_\cc,wVw^*]\in KK(\cc,\Bkom)
\]
where
\[
\begin{split}
wVw^*&=-\sigma\big(\cos(\pi q)\gradodd\big)+\sigma\big(\sin(\pi
q)\grad 1\big)V_2\\
&=-\left(\begin{matrix}0&-i\psi\big(\cos(\pi
q)\big)\\
i\psi\big(\cos(\pi q)\big)&0\end{matrix}\right)+
\left(\begin{matrix}\psi\big(\sin(\pi
q)\big)&0\\
0&\psi\big(\sin(\pi q)\big)
\end{matrix}\right)V_2\\
&=
\left(\begin{matrix}
  0 & i\psi\big(\cos(\pi q)\big)+\psi\big(\sin(\pi q)\big)(2p-1) \\
  -i\psi\big(\cos(\pi q)\big)+\psi\big(\sin(\pi q)\big)(2p-1) & 0
\end{matrix}\right)\\
&\in \ope\big((\Bkom)\oplus(\Bkom)\big)
\end{split}
\]
Here $\sigma : \ope\big((\Akom)\grad\cc_1\big)\to\ope(\Bkom)$ and
$\psi: \ope(\Akom)\to\ope(\Bkom)$ denotes the extensions as in
\cite[Proposition $2.1$]{Lance}.

Applying the isomorphism $KK(\cc,\Bkom)\cong K_0(\Bkom)$ from
Lemma \ref{isomk}
we get that the product is nothing but the element
\[
\partial\Big[\pi\big(-i\psi\big(\cos(\pi q)\big)+ \psi\big(\sin(\pi
q)\big)(2p-1)\big)\Big]\in K_0(\Bkom)
\]
Set $v=i\exp(i\pi q)=i\cos(\pi q)-\sin(\pi q)$. The element $v$ is a unitary
in $\ope(\Akom)$ and thus homotopic to $1$ so 
\[
\begin{split}
&\partial\Big[\pi\big(-i\psi\big(\cos(\pi q)\big)+ \psi\big(\sin(\pi
q)\big)(2p-1)\big)\Big]=\\
&\qquad\qquad\partial\Big[\pi\big(-i\psi\big(v\cos(\pi q)\big)+
\psi\big(v\sin(\pi q)\big)(2p-1)\big)\Big]
\end{split}
\]
Furthermore, with the same argument as in the proof of Theorem
\ref{prod}, 
we have
\[
\pi\big(\cos(\pi q)\cos(\pi q)\big)=\pi\big((1-2q)(1-2q)\big)=\pi(1)
\]
so $\cos^2(\pi q)-1\in \Akom$. Moreover $\sin(\pi q)\geq 0$ since $q=q^*$ and $\|q\|\leq 1$, 
so $\sin(\pi q)=\big(1-\cos^2(\pi q)\big)^{1/2}\in \Akom$. We thus have $v\cos(\pi q)\in (\Akom)^+$. 
By assumption $\psi(a)p-p\psi(a)\in\Bkom$ for all $a\in (\Akom)^+$ so
\begin{align*}
\pi\Big(-i\psi\big(v\cos(\pi q)\big)&+ \psi\big(v\sin(\pi
q)\big)(2p-1)\Big)\\
&=
\pi\Big(-ip\psi\big(v\cos(\pi q)\big)p-i(1-p)\psi\big(v\cos(\pi
q)\big)(1-p)\\
& +p\psi\big(v\sin(\pi q)\big)p-(1-p)\psi\big(v\sin(\pi
q)\big)(1-p)\Big)\\
&=\pi\Big(p\psi(-v^2)p+(1-p)\psi\big(v[-i\cos(\pi q)-\sin(\pi
q)]\big)(1-p) \Big)\\
&=\pi\big(p\psi(u)p+(1-p)\big)
\end{align*}
That is
\[
z=\partial\Big[\pi\big(-i\psi\big(\cos(\pi q)\big)+
\psi\big(\sin(\pi
q)\big)(2p-1)\big)\Big]\\=\partial\Big[\pi\big(p\psi(u)p+(1-p)\big)\Big]
\]
as desired.
\end{proof}


\begin{thebibliography}{99999}


\bibitem{Baaj} S. Baaj, P. Julg, 
{\em Théorie bivariante de Kasparov et opérateurs non bornés dans les 
$C^*$-modules hilbertiens},  C. R. Acad. Sci. Paris {\bf 296}, 
(1983), 139--149.
\bibitem{BCPRSW}  M-T. Benameur, A.L. Carey, J. Phillips, A. Rennie,
  F.A. Sukochev, K.P. Wojciechowski, {\em 
An Analytic Approach to Spectral Flow in von Neumann Algebras},
 in `Analysis, Geometry and Topology of
Elliptic Operators-Papers in Honour of K. P. Wojciechowski', World
Scientific, 2006,  297--352

\bibitem{Black} B. Blackadar, {\em K-Theory for Operator Algebras},  
Mathematical Sciences Research Institute
Publications, $1986$, Springer-Verlag New York.

\bibitem{BLP} B. Boo\ss-Bavnbek, M. Lesch, J. Phillips
\emph{Unbounded Fredholm operators and spectral flow}, Canad. J. Math.
{\bf 57} (2005) 225-250.
\bibitem{B1} M. Breuer, \emph{Fredholm theories in von Neumann
algebras. I}, Math. Ann., {\bfseries 178}(1968), 243-254.

\bibitem{B2} M. Breuer, \emph{Fredholm theories in von Neumann
algebras. II}, Math. Ann., {\bfseries 180}(1969), 313-325.

\bibitem{Connes} A. Connes, \emph{Noncoummutative Geometry}, Academic Press, $1994$.

\bibitem{CP1} A. Carey, J. Phillips, {\em Unbounded 
Fredholm Modules and Spectral Flow},  Can. J. Math. Vol. {\bf 50} (1998),
 673--718.

\bibitem{CP2} A.  Carey, J. Phillips,
{\em Spectral flow in $\Theta$-summable
Fredholm modules, eta invariants and the JLO cocycle}, K Theory
{\bf 31} (2004) 135-194.

\bibitem{CPRS2} A.  Carey, J. Phillips, A. Rennie, F. A. Sukochev, 
\emph{The local index formula in semifinite von Neumann algebras
I. Spectral flow}, Adv. Math, {\bf 202} No. 2 (2006)  451--516.

\bibitem{CPRS3} A. Carey, J. Phillips, A. Rennie,
  F. A. Sukochev, \emph{The local index formula in semifinite von Neumann algebras II: the even case}, 
Adv. Math, {\bf 202} No. 2 (2006)  517--554.


\bibitem{Dixmier} J. Dixmier, {\em Les algèbres d'opérateurs 
dans l'espace hilbertien (Algèbres de von Neumann)},
 Paris, Gauthiers-Villars, Editeur-Imprimeur-Libraire, 55, 
Quai Des Grands Augustins, 1957.
\bibitem{Fack} T. Fack, {\em Sur la notion de valeur caractéristique},
   J. Operator Theory, {\bf 7} (1982), 307--333.
\bibitem{Higson} N. Higson, J. Roe, {\em Analytic K-Homology}, 
Oxford University Press, (2000).
\bibitem{K} G. G. Kasparov,
{\em The operator $K$-functor and extensions of $C^*$-algebras},
Math. USSR Izv. {\bf 16} (1981), 513--572.
\bibitem{Lance} E. C. Lance, {\em Hilbert $C^*$-modules}, 
London Mathematical Society Lecture Note Series {\bf 210}, Cambridge
University Press, 1995.
\bibitem{Meise} R. Meise, D. Vogt, {\em Introduction to Functional
    Analysis},  Oxford Science Publications, Clarendon Press, 1997.

\bibitem{pr} D. Pask, A. Rennie, {\em The Noncommutative
Geometry of Graph $C^*$-Algebras I: The Index Theorem}, J. Funct. An.,
{\bf 233} (2006)  92--134
\bibitem{prs} D. Pask, A. Rennie, A. Sims, {\em The Noncommutative
Geometry of $k$-Graph $C^*$-Algebras}, math.OA/0512438 
\bibitem{Ped} G. K.
Pedersen, {\em $C^*$-Algebras and Their Automorphism Groups},  
Academic Press, London, New York, San Francisco, 1979.
\bibitem{Wegge}  N.E. Wegge-Olsen, {\em K-theory and $C^*$-algebras},
Oxford University Press, 1993.

\bibitem{Ph} J. Phillips, \emph{Self-adjoint Fredholm operators and
spectral flow}, Canad. Math. Bull., {\bf 39}(1996), 460--467.

\bibitem{Ph1} J. Phillips,
\emph{Spectral flow in type I and type II factors-a new approach},
Fields Institute Communications, {\bf vol. 17}(1997), 137--153.

\bibitem{PR}  J. Phillips, I. F. Raeburn, \emph{An index theorem for
Toeplitz operators with noncommutative symbol space}, J. Funct. Anal.,
{\bf 120} (1993) 239-263.

\bibitem{ar} A. Rennie, {\em Smoothness and Locality for 
Nonunital Spectral Triples},
$K$-Theory, {\bf 28}, (2003),  127--165

\bibitem{Roerdam} M. Rørdam, F. Larsen, N.J. Laustsen, {\em An Introduction to
$K$-Theory of $C^*$-Algebras}, Cambridge University 
Press, 2000 London Mathematical Society Student Text 49.

\bibitem{W} C. Wahl, {\em On the Noncommutative Spectral Flow},
  math.OA/0602110v3 

\end{thebibliography}
\end{document}